\documentclass{article}
\usepackage{amsmath}
\usepackage{amssymb}
\allowdisplaybreaks[4]
\usepackage{caption}
\usepackage{mathrsfs}
\usepackage{graphicx}
\usepackage{longtable,booktabs}
\usepackage{algorithm}
\usepackage{algorithmic}
\usepackage{lipsum}
\usepackage{makecell}
\usepackage{multirow}


\usepackage{setspace}
\begin{document}
	\newtheorem{Lemma}{Lemma}[section]
	\newtheorem{Proposition}{Proposition}[section]
	\newtheorem{Theorem}{Theorem}[section]
	\newtheorem{Corollary}{Corollary}[section]
	\newtheorem{Example}{Example}[section]
	\baselineskip 0.25in
	\title{\Large\bf  Measuring and aggregating $\varepsilon$-$T$-transitive fuzzy relations}
	\author{ Dechao Li$^{ \mbox{{\small a}}}$
\thanks{Email: dch1831@163. com}\qquad {Yutao Yao$^{ \mbox{{\small a}}}$} \qquad{Jingyao Duan$^{ \mbox{{\small b}}}$}\\
  $^{ \mbox{{\small a}}}${\small  School of Information and Engineering, Zhejiang Ocean University,
  Zhoushan,
  316000, China}
\\$^{ \mbox{{\small b}}}${\small Department of Mathematics and Information Sciences, BaoJi University of Arts and Sciences, Baoji, 721016, China}}
\date{}
	\maketitle
	\begin{center}
	\begin{minipage}{140mm}
		\begin{picture}(1,1)
			\line(1,0){400}
		\end{picture}

\centerline{\bf Abstract} \vskip 3mm {\qquad The transitivity of fuzzy relations plays an important role in fuzzy set theory, artificial intelligence, clustering and decision-making. However,  it is often difficult for fuzzy relations to satisfy the transitivity property in many practical applications. This has motivated researchers to investigate the degree to which a fuzzy relation is transitive. Therefore, this work first investigates two different measures of $T$-transitivity for fuzzy relations using some well-known fuzzy implications. And then, the relationship between two different degrees of transitivity is investigated. Further, the concept of an $\varepsilon$-$T$-transitive fuzzy relation is introduced, and the aggregation functions that preserve the $\varepsilon$-$T$-transitivity of fuzzy relations are characterized. Finally, the $\varepsilon$-$T$-transitive fuzzy relation is utilized to make inferences and cluster objects. Compared to finding the $T$-transitive closure, it is reasonable to cluster objects using the $\varepsilon$-$T$-transitive fuzzy relation under the permissible error.}			

\begin{picture}(1,1)
\line(1,0){400}
\end{picture}
\end{minipage}
\end{center}
\vskip 4mm
\section{Introduction}
\subsection{The motivation of this work}
\qquad As a significant foundation of fuzzy set theory,  fuzzy relations serve as a pivotal tool in various applied fields such as artificial intelligence \cite{Amini, NolaS}, approximate reasoning \cite{Castro, Zadeh1}, decision-making \cite{Fodor, Qian, Zhang}, clustering \cite{Tang, Yang}, fuzzy control \cite{ Klawonn},  image processing \cite{Bustince, Jaurrieta, Nobuhara} and data-mining \cite{Alavi}.  Particularly, these applications focus on the binary fuzzy relations defined on a universe $U$. Since a binary fuzzy relation $R$ on the universe $U$ can be regarded as a mapping $ R:U\times U\rightarrow [0,1]$ \cite{Zadeh}, we will omit the word ``binary" and only consider the binary fuzzy relation in the rest of this work. Moreover, some special properties such as reflexivity, symmetry and transitivity  are often required in the application of fuzzy relations. This motivates people to study the property of fuzzy relations\cite{Alcantara,Baets,BaetsM,Belohlavek, Dan, Demiric,Garmendia,Gupta,Hernandez,Jacas,Jacas1,Montes,Sanchez,Sun,Valverde,Wang}. Among these, the transitivity of fuzzy relation with respect to a t-norm $T$ (for short, $T$-transitivity) has been extensively  investigated.

Transitivity is essential for equivalence relations and partially ordered relations. However, in engineering and physics, people often have to deal with imprecise measurements. This imprecision can frequently be regarded as indistinguishability. If indistinguishability is modeled by an equivalence relation, transitivity will lead to a paradoxical outcome\cite{Boixader,Cock}. Consequently, multiple versions of fuzzy transitivity have been proposed\cite{Boixader,Garmendia,Klawonn}. In practical applications, ensuring a fuzzy relation satisfies transitivity is challenging. Therefore, another transitive fuzzy relation $R'$ must be found to substitute for this fuzzy relation $R$. For example, the $T$-transitive
closure $\bar{R}^T$ serves as an alternative to $R$. However, computing  $\bar{R}^T$ is computationally expensive. Its time complexity is $\bar{R}^T$ is $O(n^3)$ when $|U|=n$. And it requires storing three square $n\times n$ matrices\cite{Naessens}. Moreover, $\bar{R}^T$ does not generally ``closely" approximate $R$. To approximate better $R$, Garmendia and Recasens has proposed other transitive fuzzy relations\cite{Garmendia}.

We know that a relation $R$  is called $T$-transitive if the inequality $T(R(x,y),R(y,z))\leq R(x,z)$ holds for all $x,y,z\in U$, where $T$ is a t-norm. Obviously,  we can find some $x,y,z\in U$  such that  the inequality $T(R(x,y),R(y,z))> R(x,z)$ holds when  $R$ is not $T$-transitive. This suggests that we may gain insights into indicators of fuzzy relations linked to the properties of transitivity. Considering that the R-implication $I_T$ can be used to describe the relationship between the truth values of the antecedent and consequent in conditional statements, Wang and Xue defined a $T$-transitivity indicator using $I_T$ as follows\cite{Wang}:
\begin{equation} \alpha_{T}(R)=\underset{{x,y,z\in U}}{\inf}I_{T}(T(R(x,y),R(y,z)), R(x,z)),\end{equation}
 Later, Boixader and Recasens further studied the properties of $\alpha_T(R)$\cite{BoixaderR}. With a revised t-norm $T_c$, they demonstrated  that $\alpha_T(R)$ corresponds to the supremum of all $c\in[0,1]$ for which $R$ remains $T_c$-transitive. It should be emphasized that their analysis specifically focuses on R-implications derived from the same left-continuous t-norms as in Eq.(1).

As more fuzzy implications can be used to interpret conditional statements in fuzzy logic\cite{Baczynski}, it becomes particularly interesting to investigate the degree of $T$-transitivity for  fuzzy relations using various fuzzy implications. Therefore, this work proposes an extended concept of $T$-transitivity degree for fuzzy relations based on a generic fuzzy implication $I$, as formalized  below:
\begin{equation} \alpha_{(I,T)}(R)=\underset{{x,y,z\in U}}{\inf}I(T(R(x,y),R(y,z)), R(x,z)).\end{equation}
In this  framework, we have more choices to measure the degrees of $T$-transitivity for fuzzy relations through various fuzzy implications. This implies that this extension will be more flexible than Eq.(1). Clearly, $\alpha_{(I,T)}(R)$  constitutes an extension of $\alpha_T(R)$ when $I=I_T$. Thus, we say that $R$ is $\varepsilon$-$T$-transitive if $\alpha_{(I,T)}(R)$ is greater than or equal to a threshold $\varepsilon$ in this work. And then we will obtain the concept of $\varepsilon$-fuzzy equivalence relations if both the degree of symmetry and the degree of $T$-transitivity are greater than or equal to a given threshold $\varepsilon$.

In addition, the fuzzy implication has been  employed to measure the similarity of two fuzzy sets\cite{Bandler}. Specifically, the similarity measure $S_I$ is defined as $S_I(A,B)=
\mathop{\inf}\limits_{x\in U}\min(I(A(x),B(x)),$ $I(B(x),A(x)))$, where $I$ denotes a fuzzy implication. Since the fuzzy relation $R$ can be viewed as a fuzzy set on $U^2$, the similarity  of two fuzzy relations $R$ and $R'$ on the\vspace{1mm} same universe $U$ can be analogously defined as $S_I(R,R')=\mathop{\inf}\limits_{x,y\in U}\min (I(R(x,y),R'(x,y)), I(R'(x,y),R(x,y)))$. As\vspace{1mm} revealed by Boixader and Recasens\cite{BoixaderR}, $\alpha_T(R)$ ia related to a $T_c$-transitive fuzzy relation. This motivates employing an alternative transitive relation to characterize the transitivity degree for fuzzy relations. Therefore, we can also define alternative measure of $T$-transitivity of fuzzy relations using a transitive fuzzy relation.  Let $\mathcal{R_T}$ denote the family of all $T$-transitive fuzzy relations on $U$. Another degree of $T$-transitivity for fuzzy relation $R$ is defined as follows:
\begin{equation} \alpha_{S_I}(R)=\sup_{R'\in\mathcal{R_T}}S_I(R,R').\end{equation}
This naturally leads to an interesting topic to investigate the fundamental connections between $\alpha_{(I,T)}(R)$ and $\alpha_{S_I}(R)$.

How to aggregate some incoming data into a single output is a crucial issue in information fusion. As fundamental tools, aggregation functions have been widely utilized in many intelligent systems. Information aggregation often requires the output to share the same properties as the inputs. Therefore, it becomes urgent to select appropriate aggregation operators for combining fuzzy relations in many practical applications. This motivates our insight into the preservation of some properties during the process aggregating fuzzy relations. Particularly, we rigorously consider the conditions under which $\varepsilon$-$T$-transitivity is preserved when aggregating a series of $\varepsilon$-$T$-transitive fuzzy relations.
\subsection{Contribution of this work}
\qquad Building upon the aforementioned arguments, we can effectively analyze the degree of transitivity with some fuzzy implications. To advance the capability of fuzzy relation in artificial intelligence, clustering and decision-making, this work will methodically examine the two distinct measures of $T$-transitivity employing various fuzzy implications. Subsequently, we will investigate the practical implementations of applications of $\varepsilon$-$T$-transitive fuzzy relation in fuzzy reasoning and  clustering. The work mainly contributes to

(1) Giving two  distinct degrees of $T$-transitivity using various fuzzy implications.

(2) Aggregating the $\varepsilon$-$T$-transitive fuzzy relations.

(3) Applying the $\varepsilon$-$T$-transitive fuzzy relation to data clustering.

Thus, we organize this work as follows. In Section 2, some basic notions and properties of aggregation functions, fuzzy implications and
fuzzy relations are recalled.  Section 3 introduces and analyzes two distinct measures of $T$-transitivity for fuzzy relations using various fuzzy implications, along with their interrelationships. Section 4 examines the aggregation properties of $\varepsilon$-$T$-transitive fuzzy relations. Section 5 demonstrates practical applications of $\varepsilon$-$T$-transitive fuzzy relation in approximate reasoning and clustering.
	\section{Preliminary}
	\qquad To better present the content of this work, this section will recall some concepts and properties used in  this work.
\\{\bf Definition 2.1}\cite{Grabisch}  $G_n:[0,1]^n\rightarrow [0,1]$ is an $n$-ary aggregation function if
	
	(G1)  $G_n$ is non-decreasing in every variable,
	
	(G2) $G_n(0, \cdots, 0) =0$ and  $G_n(1,\cdots, 1)=1$.
\\{\bf Example 2.2}\cite{Grabisch} The following are some basic $n$-ary aggregation functions:
\begin{itemize}
  \item The minimum function ${\rm Min}(x_1,\cdots,x_n)=\min\{x_1,\cdots,x_n\}$;
  \item The maximum function ${\rm Max}(x_1,\cdots,x_n)=\max\{x_1,\cdots,x_n\}$;
  \item Triangular norms and conorms (see Remark 1);
  \item Weighted quasi-arithmetic mean (WQAM) $G_{f,\omega}(x_1,\cdots,x_n)=f^{-1}\left(\mathop{\sum}\limits_{i=1}^n\omega_if(x_i)\right)$, where $f:[0, 1]\rightarrow$\vspace{1mm} $(-\infty, +\infty)$ is a strictly
monotone continuous mapping.
\end{itemize}
{\bf Definition 2.3}\cite{Klement} A t-norm $T$(t-conorm $S$) is an associative and commutative binary aggregation function having a neutral element 1(0).
\\{\bf Remark 1.} The associativity of t-norms (t-conorms) allows their extension to $n$-ary aggregation functions, which are respectively defined as
$T_n(x_1,\cdots,x_n)=T(T(x_1,\cdots, x_{n-1}),x_n)$ and $S_n(x_1,\cdots,x_n)=S(S(x_1,\cdots, x_{n-1}),x_n)$.
\\{\bf Example 2.4}\cite{Klement} The four well-known t-norms(t-conorms) are shown in the following\vspace{2mm}
\begin{itemize}	
\item $T_{\rm{D}}(x,y)=\left\{\begin{array}{ll}
		0 & x,y \in [0,1)^2\vspace{1mm}\\
		\min(x,y) & \textmd{otherwise}
	\end{array}\right.$ and
	$S_{\rm{D}}(x,y)=\left\{\begin{array}{ll}
		1 & x,y \in (0,1]^2\vspace{1mm}\\
		\max(x,y) & \textmd{ otherwise}
	\end{array}\right.$;

\item $T_{\rm{L}}(x,y)=\max(x+y-1,0)$ and $S_{\rm{L}}(x,y)=\min(x+y,1)$;

\item $T_{\rm{M}}(x,y)=\min(x,y)$ and
	$S_{\rm{M}}(x,y)=\max(x,y)$;
	
\item $T_{\rm{P}}(x,y)=xy$ and $S_{\rm{P}}(x,y)=x+y-xy$.
\end{itemize}
{\bf Remark 1.} The associativity of t-norm (t-conorm) allows us  respectively defined a $n$-ary operation as
$T_n(x_1,\cdots,x_n)=T(T(x_1,\cdots, x_{n-1}),x_n)$ and $S_n(x_1,\cdots,x_n)=S(S(x_1,\cdots, x_{n-1}),x_n)$.
\\{\bf Definition 2.5}\cite{Klement} i.  $x\in (0,1)$ is a nilpotent element of $T$ if  $T_n(x,\cdots,x)=0$ holds for an $n\in \mathbb{N}$.\vspace{1mm}

ii. A continuous $T$ is nilpotent if  every $x\in (0,1)$ is its nilpotent element;

iii. $T$ is Archimedean if  we can find an $n\in \mathbb{N}$ such that the inequality $T_n(x,\cdots,x)<y$ holds for all $x$ and $y\in (0,1)$;\vspace{1mm}

iv. $T$ is strict if it is strictly monotone (that is, $T(x,y)<T(x,z)$ holds for any $y<z$ and $x>0$) and continuous.
\\{\bf Proposition 2.6}\cite{Klement} i. A continuous $T$ is Archimedean $\Longleftrightarrow T(x,y)=t^{-1}(\min(t(0),t(x)+t(y)))$ with an additive generator $t$, where $t$ is a continuous and strictly decreasing mapping from [0,1] to $[0,+\infty]$ with $t(1)=0$.

ii. $T$ is strict $\Longleftrightarrow t(0)=\infty\Longleftrightarrow$  $T=(T_{\rm P})_\varphi$, where $(T_{\rm P})_\varphi(x,y)=\varphi^{-1}(\varphi(x)\varphi(y))$ and $\varphi$ is a strictly increasing bijection on [0,1];

iii. $T$ is nilpotent $\Longleftrightarrow t(0)<\infty\Longleftrightarrow$  $T=(T_{\rm L})_\varphi$.
\\{\bf Definition 2.7}\cite{Zadeh} A fuzzy relation from $U$ to $V$ is a fuzzy set $R: U\times V\rightarrow [0,1]$. Especially, $R$ is a fuzzy relation on $U$ if $U=V$.
\\{\bf Definition 2.8}\cite{Cock,Zadeh}  We say $R$ is

i. reflexive if $R(x,x)=1$ holds for any $x\in U$;

ii.  symmetric if $R(x,y)=R(y,x)$ holds for any $x, y\in U$;

iii. $T$-transitive if the inequality $R(x,z)\geq T(R(x,y),R(y,z))$ holds for any $x,y,z\in U$.

Further, a reflexive, symmetric and  $T$-transitive fuzzy relation $R$ is  $T$-equivalent.
\\{\bf Remark 2.} The $T$-equivalent fuzzy relation is also called as a similarity fuzzy relation or $T$-indistinguishability fuzzy relation by some researchers\cite{Valverde, Zadeh}.
\\{\bf Definition 2.9}\cite{Baczynski} A fuzzy implication $I$ is a binary operation on [0,1] satisfied the following

(I1) $I(x,z)\geq I(y, z)$ if $x\leq y$;

(I2) $I(x,y)\leq I(x, z)$ if $y\leq z$;

(I3) $I(0,0)=1$;

(I4) $I(1,1)=1$;

(I5) $I(1,0)=0$.
\\{\bf Definition 2.10} \cite{Baczynski}  We say that $I$ fulfills

(OP) For all $x,y\in [0,1]$, $I(x, y)=1\Longleftrightarrow x\leq y$.
\\{\bf Definition 2.11} \cite{Baczynski,Fodor}  An $(S,N)$-implication is a mapping $I_{(S,N)}:[0,1]^2\rightarrow[0,1]$  defined as $I_{(S,N)}(x,y)=S(N(x), y)$.
\\{\bf Definition 2.12}\cite{Baczynski} A R-implication $I_T$ generated by $T$ is $I_T(x,y) = \sup\{ z \in [0, 1]\mid T(x, z) \leq y\}$.
\\{\bf Proposition 2.13} \cite{Baczynski} Let $T$ be left-continuous. We have

i. $y\leq I_T(x,z)\Longleftrightarrow T(x,y)\leq z$;

ii. $I_T$ satisfies OP;\vspace{1mm}

iii. $I_T\left(\mathop{\sup}\limits_{\lambda\in \Lambda}x_\lambda,y\right)=\mathop{\inf}\limits_{\lambda\in \Lambda}I_T(x_\lambda,y)$ and $I_T\left(x,\mathop{\inf}\limits_{\lambda\in \Lambda}y_\lambda\right)=\mathop{\inf}\limits_{\lambda\in \Lambda}I_T(x,y_\lambda)$.\vspace{1mm}
\\{\bf Definition 2.14}\cite{Baczynski} A QL-operation $I_{(T,S,N)}$  is defined as $I_{(N, T,S)}(x, y) =S(N(x), T(x, y))$. Especially, it is  a QL-implication if it fulfills I1-I5.
\\{\bf Definition 2.15}\cite{Massanet} $I^T(x, y) = \sup\{r\in [0, 1] | y^{(r)}
_T \geq x\}$  is  a $T$-power implication.
\\{\bf Proposition 2.16}\cite{Massanet} Let $T$ be Archimedean with additive generator $t$. Then
 $I^T(x, y)=$\vspace{1mm} $\left\{\begin{array}{ll}
                                                                                                 1 & x\leq y\\
                                                                                                 \frac{t(x)}{t(y)}& x>y
                                                                                               \end{array}\right..$\vspace{1mm}
\\{\bf Definition 2.17}\cite{Yager} Suppose that $g:[0,1]\rightarrow [0,+\infty]$ is a strict increasing and continuous mapping with\vspace{1ex}
$g(0)=0$. A $g$-implication with a $g$-generator is defined by
$I_{g}(x,y)=g^{(-1)}\left(\frac{g(y)}{x}\right)$ with $0\times \infty=\infty$, where $g^{(-1)}(x)=\left\{
\begin {array}{ll}
 g^{(-1)}(x) &  x\leq g(1)\vspace{1mm}\\
  1 &  \textmd{otherwise}
\end{array}\right.$.\vspace{1mm}
\\{\bf Definition 2.18}\cite{Grzegorzewski} Let  $C$ be a copula. A probabilistic implication is\vspace{2mm} $I_C(x,y)=\left\{\begin{array}{ll}
                                                                                   \frac{C(x,y)}{x} & x>0 \vspace{1mm}\\
                                                                                    1 & x=0
                                                                                  \end{array}\right.$
if it satisfies I1.\vspace{1mm}
\\{\bf Definition 2.19}\cite{Grzegorzewski} A probabilistic S-implication is
$\tilde{I}_C(x,y)=C(x,y)-x+1$.
  \section{The measures of $T$-transitivity for fuzzy relation}
   \qquad In this section, we systematically analyze two distinct measures of $T$-transitivity for fuzzy relations using various fuzzy implications introduced in Section 2.
\subsection{Fuzzy implications-based $T$-transitivity measurement via $\alpha_{(I,T)}(R)$}
\qquad Let us first formalize the concept of $\varepsilon$-$T$-transitive fuzzy relation through the $\alpha_{(I,T)}(R)$ measure as follows.
\\{\bf Definition 3.1} Let $\varepsilon\in [0,1]$ be arbitrary a fixed constant. We say that  $R$ is $\varepsilon$-$T$-transitive if $\alpha_{(I,T)}(R)\geq \varepsilon$.
\\{\bf Remark 3.} i. According to Definition 3.1, $\alpha_{(I,T)}(R)\geq \varepsilon=1\Longleftrightarrow R$ is $T$-transitive. Therefore, we will omit the prefix ``$\varepsilon$-"
in Definition 3.1 when $\varepsilon=1$.

ii. When the fuzzy implication $I$ in Eq.(2) satisfies OP, the following statement holds: $R$ is $T$-transitive $\Longleftrightarrow\alpha_{(I,T)}(R)=1$. However,  $\alpha_I(R)$ may deviate from 1 for  a $T$-transitive fuzzy relation $R$ when $I$ violates OP. Consequently, we explicitly require  $I$ in Eq.(2) to satisfy OP throughout this work.

iii. As mentioned above, the $T$-transitivity of $R$ brings to light the relationship between $R(x,z)$ and $T(R(x,y),R(y,z))$.  Similarly, we examine the structural properties of $\varepsilon$-$T$-transitive fuzzy relation though investigating the relationship between $R(x,z)$ and $T(R(x,y),R(y,z))$ for various fuzzy implications introduced in Section 2.
  \\{\bf Proposition 3.2} Let $I$ be a continuous $(S,N)$-implication in Eq.(2). Then, $R$ is $\varepsilon$-$T$-transitive $\Longleftrightarrow T(R(x,y),R(y,z))\leq \varphi^{-1}(\min (1,1+\varphi(R(x,z))-\varphi(\varepsilon)))$ holds for any $x,y,z\in U$.
  \\{\bf Proof.} In this case, $I$ can be rewritten as $I_{(S,N)}(x,y)=\varphi^{-1}(\min(1,1-\varphi(x)+\varphi(y)))$\cite{Baczynski,Fodor}, where $\varphi$ is an increasing bijection on [0,1].

  ($\Longrightarrow$) Suppose that $R$ is  $\varepsilon$-$T$-transitive. We have $I_{(S,N)}(T(R(x,y),R(y,z)),R(x,z))=\varphi^{-1}(\min(1,1-\varphi(T(R(x,y),R(y,z)))+\varphi(R(x,z))))\geq \varepsilon$ for all $x,y,z\in U$. This implies that $1-\varphi(T(R(x,y),R(y,z)))+\varphi(R(x,z))\geq\varphi(\varepsilon)$ holds. Further, this inequality can be rewritten as $\varphi(T(R(x,y),R(y,z)))\leq1-\varphi(\varepsilon)+\varphi(R(x,z))$. We therefore have $T(R(x,y),R(y,z))\leq \varphi^{-1}(\min (1,1+\varphi(R(x,z))-\varphi(\varepsilon)))$.

  ($\Longleftarrow$) Let $T(R(x,y),R(y,z))\leq \varphi^{-1}(\min (1,1+\varphi(R(x,z))-\varphi(\varepsilon)))$ for all $x,y,z\in U$. This means that $\varphi(T(R(x,y),R(y,z)))\leq 1+\varphi(R(x,z))-\varphi(\varepsilon)$ holds. We further have  $1-\varphi(T(R(x,y),R(y,z)))+\varphi(R(x,z))\geq\varphi(\varepsilon)$. Therefore, $\min(1,1-\varphi(T(R(x,y),R(y,z)))+\varphi(R(x,z)))\geq\varphi(\varepsilon)$ holds. Thus, we obtain $I_{(S,N)}(T(R(x,y),R(y,z)),R(x,z))\geq \varepsilon$. That is, $R$ is $\varepsilon$-$T$-transitive.$\hfill\square$
  \\{\bf Remark 4.} i.  For fixed $x, z\in U$, consider the function $f(\varepsilon)=\varphi^{-1}(\min (1,1+\varphi(R(x,z))-\varphi(\varepsilon)))$. Obviously, $f(\varepsilon)$ is continuous decreasing on [0,1] and $\mathop{\lim}\limits_{\varepsilon\rightarrow 1^-}f(\varepsilon)=R(x,z)$. Consequently, we observe that $R$ attains higher $T$-transitivity when $\alpha_{(I,T)}(R)$ increases. This implies that $R$ approaches maximal $T$-transitivity as $\varepsilon\rightarrow 1^-$, which aligns precisely with Definition 3.1.

  ii. Let $U_0^2=\{(x,z)\in U^2|R(x,z)<T(R(x,y),R(y,z))\leq f(\varepsilon)\ {\rm for\ some}\ y\in U\}\ (\subseteq U^2)$. As demonstrated in Proposition 3.2, the cardinality of $U^2_0$ decreases monotonically as $\varepsilon\rightarrow 1^-$.  Notably, $U_0^2=\emptyset\Longleftrightarrow R$ is $T$-transitive.

  Further, let $T=(T_{\rm L})_\varphi$. The original inequality  $T(R(x,y),R(y,z))\leq \varphi^{-1}(\min (1,1+\varphi(R(x,z))-\varphi(\varepsilon)))$ can be reformulated as $\max(0,\varphi(R(x,y))+\varphi(R(y,z))-1)\leq \min (1,1+\varphi(R(x,z))-\varphi(\varepsilon))$. This inequality automatically holds when $\varphi(R(x,y))+\varphi(R(y,z))\leq1$. In the alternative case, it cab be reduced as $\varphi(R(x,y))+\varphi(R(y,z))-1\leq 1+\varphi(R(x,z))-\varphi(\varepsilon)$. Further, we have $(\varphi(R(x,y))+\varphi(\varepsilon)-1)+(\varphi(R(y,z))+\varphi(\varepsilon)-1)-1\leq \varphi(R(x,z))+\varphi(\varepsilon)-1$. Defining $R_{(({T_{\rm L}})_\varphi,\varepsilon)}(x,y)=(T_{\rm L})_\varphi(R(x,y),\varepsilon)$,  we formally state the consequent corollary.
  \\{\bf Corollary 3.3}  Let $I$ be a continuous $(S,N)$-implication in Eq.(2). $R_{(({T_{\rm L}})_\varphi,\varepsilon)}$ is $(T_{\rm L})_\varphi$-transitive if $R$ is a $\varepsilon$-$(T_{\rm L})_\varphi$-transitive fuzzy relation.
  \\{\bf Remark 5.} i. However, this condition is not necessary as shown in Example 3.4.

  ii. The condition $R(x,y)\geq \varphi^{-1}(1-\varphi(\varepsilon))$ provides a necessary condition as shown in Example 3.5.
  \\{\bf Example 3.4} Considering the fuzzy relation $R=\left(\begin{array}{cc}
                                                                0.2 & 0.8 \\
                                                                0.8& 0.2
                                                              \end{array}\right)$. Let $\varepsilon=0.7$ and $\varphi(x)=x$.\vspace{1mm} We have $R_{(\varphi,T_{\rm L},\varepsilon)}=\left(\begin{array}{cc}
                                                                0 & 0.5 \\
                                                                0.5& 0
                                                              \end{array}\right)$. Obviously, $R_{(\varphi,T_{\rm L},\varepsilon)}$ is $T_{\rm L}$-transitive. However, $\alpha_{(I,T)}(R)=0.6$\vspace{1mm} implies that $R$ is not 0.7-$T_{\rm L}$-transitive.
   \\{\bf Example 3.5} Considering the fuzzy relation $R=\left(\begin{array}{cc}
                                                                0.2 & 0.8 \\
                                                                0.8& 0.2
                                                              \end{array}\right)$ again. Let $\varphi(x)=x^2$. In\vspace{1mm} order to ensure that  $R(x,y)\geq \varphi^{-1}(1-\varphi(\varepsilon))$ holds for all $x,y\in [0,1]$, take $\varepsilon=\sqrt{1-0.8^2}$\vspace{1mm} $=0.6$. We have $R_{(({T_{\rm L}})_\varphi,\varepsilon)}=\left(\begin{array}{cc}
                                                                0 & 0.28 \\
                                                                0.28& 0
                                                              \end{array}\right)$. Obviously, $R_{(({T_{\rm L}})_\varphi,\varepsilon)}$ is $(T_{\rm L})\varphi$-transitive. In this\vspace{1mm} case, $\alpha_{(I,T)}(R)=0.6$. That is, $R$ is  $0.6$-$(T_{\rm L})_\varphi$-transitive.
  \\{\bf Proposition 3.6} Let $I$ be a R-implication generated by a left-continuous $T'$ in Eq.(2). Then, $R$ is $\varepsilon$-$T$-transitive $\Longleftrightarrow T'(\varepsilon, T(R(x,y),R(y,z)))\leq R(x,z)$ holds for all $x,y,z\in U$.
  \\{\bf Proof.} $(\Longrightarrow)$ Let $R$ be a $\varepsilon$-$T$-transitive fuzzy relation. We have $\alpha_{(I,T)}(R)=\mathop{\inf}\limits_{x,y,z\in U}I_T(T(R(x,y),$\vspace{1mm} $R(y,z)),R(x,z))\geq \varepsilon$. The left continuity of $T'$ implies that  $I_T(T(R(x,y),R(y,z)),R(x,z))\geq \varepsilon\Longleftrightarrow T'(\varepsilon, T(R(x,$ $y),R(y,z)))\leq R(x,z)$.

  $(\Longleftarrow)$ For any $x,y,z\in U$, suppose that $T'(\varepsilon, T(R(x,y),R(y,z)))\leq R(x,z)$ holds. Since $T'$ is left-continuous, we have $I_T(T(R(x,y),R(y,z)),R(x,z))\geq \varepsilon$. Thus, $\alpha_{(I,T)}(R)=\mathop{\inf}\limits_{x,y,z\in U}I_T(T(R(x,y),R(y,z)),R(x,z))\geq \varepsilon$ holds. This means that $R$ is $\varepsilon$-$T$-transitive.$\hfill\square$\vspace{1mm}

  As $T'$ is a left-continuous t-norm, we have $\mathop{\sup}\limits_{y\in U}T'(\varepsilon, T(R(x,y),R(y,z)))=T'(\varepsilon, \mathop{\sup}\limits_{y\in U}T(R(x,y),$\vspace{1mm} $R(y,z)))$. Let $R_T(x,y)=\mathop{\sup}\limits_{z\in U}T(R(x,z),R(z,y))$. Defining $R_{(T,T',\varepsilon)}(x,y)=T'(R_T(x,y),\varepsilon)$, we\vspace{1mm} have the fact that $R$ is $\varepsilon$-$T$-transitive  $\Longleftrightarrow R_{(T,T',\varepsilon)}\leq R$. Moreover, $\mathop{\lim}\limits_{\varepsilon\rightarrow 1^-}R_{(T,T',\varepsilon)}(x,y)=T'(R_T(x,y),\mathop{\lim}\limits_{\varepsilon\rightarrow 1^-}\varepsilon)=R_T(x,y)$ holds. This implies that the the distance between $R_{(T,T',\varepsilon)}$ and $R_T$ decreases as $\varepsilon$   approaches 1. Especially, $R_T$ serves as the reachable least upper bound of the family $\{R_{(T,T',\varepsilon)}\}_{\varepsilon\in[0,1]}$. We formally state the following result.
  \\{\bf Corollary 3.7} Let $I$ be a R-implication generated by a continuous Archimedean $T'$ in Eq.(2). $T(R(x,y),R(y,z))\leq t^{-1}(t(\varepsilon)+t(R(x,z)))$ holds for any $x,y,z\in U$ if $R$ is $\varepsilon$-$T$-transitive, where $t$ is the additive generator of $T$.
  \\{\bf Proof.} In this case, $I_T(x,y)=t^{-1}(\max(0,t(y)-t(x)))$. Since $R$ is  $\varepsilon$-$T$-transitive, $\alpha_{(I,T)}(R)\geq \varepsilon$ implies that $t^{-1}(\max(0,t(T(R(x,y),R(y,z)))-t(R(x,z))))\geq \varepsilon$ holds. Thus, we have $T(R(x,y),$ $R(y,z))\leq t^{-1}(t(\varepsilon)+t(R(x,z)))$ for all $x,y,z\in U$.$\hfill\square$

 Especially, the following result holds if $I_T$ is generated by the same t-norm $T$.
  \\{\bf Corollary 3.8} Let $I$ be a R-implication generated by a left-continuous $T$ in Eq.(2). $T(T(R(x,$ $y),\varepsilon), T(R(y,z),\varepsilon))\leq T(R(x,z),\varepsilon)$ holds for any $x,y,z\in U$ if $R$ is $\varepsilon$-$T$-transitive.
  \\{\bf Proof.} By Proposition 3.6, we have $T(\varepsilon, T(R(x,y),R(y,z)))\leq R(x,z)$. This implies that $T(T(R(x,y),\varepsilon), T(R(y,z),\varepsilon))\leq T(R(x,z),\varepsilon)$ holds for any $x,y,z\in U$.$\hfill\square$
  \\{\bf Remark 6.} i. Let $R_{(T,\varepsilon)} (x,y)=T(R(x,y),\varepsilon)$. Corollary 3.8 implies that $R_{(T,\varepsilon)}$ is $T$-transitive when $R$ is a $\varepsilon$-$T$-transitive fuzzy relation.

  ii. Similar to Corollary 3.3,  this condition is not necessary.

  Let $I_{(N,T,S)}$ denote a  QL-operation generated by a strict $N$, a continuous  $T$ and a continuous  $S$. The necessary condition for $I_{(N,T,S)}$ to qualify as a fuzzy implication requires $S=(S_{\rm L})_\varphi$ and $N(x)\geq \varphi^{-1}(1-\varphi(x))$\cite{Baczynski,Shi}. Consequently, we focus on the QL-implication $I_{(N,T,S)}(x,y)=(S_{\rm L})_\varphi(N(x),T(x,y))$, where $N(x)\geq \varphi^{-1}(1-\varphi(x))$. When $I_{(N,T,S)}$ fulfills OP,  further analysis reveals that $T$ is neither e  nilpotent t-norm nor a collection of unique continuous
Archimedean t-norms\cite{Shi}. Therefore, we restrict our analysis to the two specific t-norms: $T=T_{\rm M}$ and $T=T_{\rm P}$. Under these constraints, the following statements can be obtained.
 \\{\bf Proposition 3.9} Let $I$ be specified as the QL-implication $I_{(N,T_{\rm M},S)}$ in Eq.(2). Then, $R$ is $\varepsilon$-$T$-transitive $\Longleftrightarrow T(R(x,$ $y),R(y,z))\leq N^{-1}(\varphi^{-1}(\max(0,\varphi(\varepsilon)-\varphi(R(x,z)))))$ holds for all $(x,z)\in U_1^2$ and $y\in U$, where $U_1^2=\{(x,z)\in U^2|R(x,z)<T(R(x,y),R(y,z))\ {\rm for\ some} \ y\in U\}$.
  \\{\bf Proof.}  $\alpha_{(I,T)}(R)\geq \varepsilon \Longleftrightarrow \varphi^{-1}(\min(1,\varphi(N(T(R(x,y),R(y,z))))+\varphi(\min(T(R(x,y),R(y,z)),$ $R(x,z)))))\geq \varepsilon$ holds for all $x,y,z\in U$ $\Longleftrightarrow \varphi^{-1}(\min(1,\varphi(N(T(R(x,y),R(y,z))))+\varphi(R(x,z))))$ $\geq \varepsilon$  for all $(x,z)\in U_1^2$ and $y\in U$$\Longleftrightarrow T(R(x,y),$ $R(y,z))\leq N^{-1}(\varphi^{-1}(\max(0,\varphi(\varepsilon)-\varphi(R(x,z)))))$ holds for all $(x,z)\in U_1^2$ and $y\in U$.$\hfill\square$

For others fuzzy implications, analogous conditions for $\varepsilon$-$T$-transitivity of $R$ can be derived. Consequently, these conditions are summarized in the subsequent tables.
 \begin{center}
  \mbox{\bf{\small Table 1 The necessary and sufficient conditions for $\varepsilon$-$T$-transitivity of $R$ }}\vspace{1mm}
 \end{center}
 \begin{center}
\begin{tabular}{cc}
    \toprule[1pt]
    Fuzzy implications&  Necessary and sufficient conditions\\
    \midrule
$I_{(N,T_{\rm P},S)}$ & $T(R(x,y),R(y,z))\leq \min\left(1,\frac{R(x,z)+\sqrt{R^2(x,z)+4(1-\varepsilon)}}{2}\right)$\vspace{1mm}\\
   $I^T$ with an additive generator $t$ & $T(R(x,y),R(y,z))\leq t^{-1}(\varepsilon t(R(x,z)))$\vspace{1mm}\\
   $I_g$ & $T(R(x,y),R(y,z))\leq\min\left(1,\frac{R(x,z)}{\varepsilon}\right)$\vspace{1mm}\\
   $I_C$ & $T(R(x,y),R(y,z))\leq \min\left(1,\frac{R(x,z)}{\varepsilon}\right)$\vspace{1mm}\\
    $\widetilde{I}_C$& $T(R(x,y),R(y,z))\leq \min(1,1-\varepsilon+R(x,z))$ \vspace{1mm}\\
    \bottomrule[1pt]
  \end{tabular}
\end{center}\vspace{2mm}
 \begin{center}
  \mbox{\bf{\small Table 2 The results when $R$ is $\varepsilon$-$T$-transitive}}\vspace{1mm}
 \end{center}
 \begin{center}
\begin{tabular}{ccc}
    \toprule[1pt]
    Fuzzy implications&  $\varepsilon$-$T$-transitivity of $R$ &Results\\
    \midrule
  $I_{(N,T_{\rm M},S)}$&$\varepsilon$-$(T_{\rm L})_\varphi$-transitive &$R_{((T_{\rm L})_\varphi,\varepsilon)}$ is $(T_{\rm L})_\varphi$-transitive\vspace{1mm}\\
   $I_g$ & $\varepsilon$-$T_{\rm P}$-transitive&$R_{(T_{\rm P},\varepsilon)}$ is $T_{\rm P}$-transitive\vspace{1mm}\\
   $I_C$ & $\varepsilon$-$T_{\rm P}$-transitive& $R_{(T_{\rm P},\varepsilon)}$ is $T_{\rm P}$-transitive\vspace{1mm}\\
    $\widetilde{I}_C$& $\varepsilon$-$T_{\rm L}$-transitive& $R_{(T_{\rm L},\varepsilon)}$ is $T_{\rm L}$-transitive\vspace{1mm}\\
    \bottomrule[1pt]
  \end{tabular}
\end{center}\vspace{1mm}
 \subsection{ The degree of $T$-transitivity quantified by $\alpha_{S_I}(R)$}
  \qquad This subsection investigate the properties of $\alpha_{S_I}(R)$ and  its relationship with $\alpha_{(I,T)}(R)$. To ensure consistency, we  require that $I$ in Eq.(3) satisfies OP.  Under this constraint, $\alpha_{S_I}(R)$ can be\vspace{1mm} reformulated as $\alpha_{S_I}(R)=\mathop{\sup}\limits_{R'\in\mathcal{R_T}}\left(\min\left(\mathop{\inf}\limits_{R>R'}I(R(x,y),R'(x,y)),\mathop{\inf}\limits_{R\leq R'}I(R'(x,y),R(x,y))\right)\right)$.\vspace{1.5mm} The following result is then derived.
   \\{\bf Proposition 3.10} Let $I$ be a continuous fuzzy implication in Eq.(3). Then, $\alpha_{S_I}(R)=1\Longleftrightarrow R$ is $T$-transitive.
   \\{\bf Proof.} ($\Longleftarrow$) Obviously.

   ($\Longrightarrow$) Let $\alpha_{S_I}(R)=1$. Then, we have $1=\alpha_{S_I}(R)\leq\mathop{\sup}\limits_{\delta\in[0,1]}\mathop{\inf}\limits_{x,z\in U}\min(I(R(x,z),\delta), I(\delta,R(x,$ $z)))=\mathop{\sup}\limits_{\delta\in[0,1]}\min\left(\mathop{\inf}\limits_{R(x,z)\leq T(R(x,y),R(y,z))}\min(I(R(x,z),\delta), I(\delta,R(x,z))),\right.\mathop{\inf}\limits_{R(x,z)> T(R(x,y),R(y,z))}$\vspace{1mm} $\min(I(R(x,z),\delta), I(\delta,R(x,z))))\leq \mathop{\sup}\limits_{\delta\in[0,1]}\min\left(\mathop{\inf}\limits_{R(x,z)\leq T(R(x,y),R(y,z))}\min(I(R(x,z),\delta), I(\delta,T(R(\right.$\vspace{1mm} $\left.x,y),R(y,z)))),\mathop{\inf}\limits_{R(x,z)> T(R(x,y),R(y,z))}\min(I(T(R(x,y),R(y,z)),\delta), I(\delta,R(x,z)))\right)\leq\min\left(\mathop{\sup}\limits_{\delta\in[0,1]}\right.$\vspace{1.5mm} $\mathop{\inf}\limits_{R(x,z)\leq T(R(x,y),R(y,z))}\min(I(R(x,z),\delta), I(\delta,T(R(x,y),R(y,z))),\ \mathop{\sup}\limits_{\delta\in[0,1]}\mathop{\inf}\limits_{R(x,z)> T(R(x,y),R(y,z))}\min$\vspace{1mm} $(I(T(R(x,y),R(y,z)),\delta), I(\delta,R(x,z))))=\mathop{\sup}\limits_{\delta\in[0,1]}\mathop{\inf}\limits_{R(x,z)> T(R(x,y),R(y,z))}\min(I(T(R(x,y),R(y,z)),$\vspace{1mm} $\delta), I(\delta,R(x,z)))$. We can assert that $R$ is
    $T$-transitive.  Otherwise, the inequality $T(R(x_0, y_0),$ $R(y_0, z_0))> R(x_0, z_0)$ holds for some $x_0, y_0, z_0 \in U$. Further, we define a mapping $f$ on [0,1] as
  $$f(\delta)=\min(I(T(R(x_0, y_0), R(y_0, z_0)),\delta),I(\delta,R(x_0,z_0))).$$
 It is evident that $f$ is continuous. Let $\Delta=\{\delta\in$\vspace{1mm} $[0,1]|I(T(R(x_0, y_0), R(y_0, z_0)),\delta)=I(\delta,R(x_0,$ $z_0))\}$. Since $I$ is continuous, we have $\mathop{\sup}\limits_{\delta\in[0,1]}f(\delta)=$\vspace{1mm} $\mathop{\sup}\limits_{\delta\in \Delta} f(\delta)\leq I(T(R(x_0, y_0),R(y_0, z_0)),R(x_0,$\vspace{1mm} $z_0))<1$. This implies that $\alpha_{S_I}(R)\leq\mathop{\sup}\limits_{\delta\in[0,1]}f(\delta)<1$, which contradicts $\alpha_{S_I}(R)=1$. Therefore, $R$ is\vspace{1mm}  $T$-transitive. $\hfill\square$

  According to the proof of Proposition 3.10, the inequality  $\alpha_{S_I}(R)\leq\mathop{\sup}\limits_{\delta\in[0,1]}f(\delta)<I(T(R(x,$\vspace{1mm} $y),R(y, z)),R(x,z))$ holds for any $x,y,z\in U$. This leads to the following result.
  \\{\bf Proposition 3.11}  Let $I$ be a continuous fuzzy implication. Then, we have $\alpha_{S_I}(R) \leq \alpha_{(I,T)}(R)$.
 \\{\bf Proof.} We have $\alpha_{S_I}(R)\leq\mathop{\sup}\limits_{\delta\in[0,1]}f(\delta)<I(T(R(x, y),R(y, z)),R(x,z))\leq \mathop{\inf}\limits_{x,y,z\in U}I(T(R(x, y),$\vspace{1mm} $R(y, z)),R(x,z))=\alpha_{(I,T)}(R)$.  $\hfill\square$

  To further investigate the relationship between $\alpha_{S_I}(R)$ and $\alpha_{(I,T)}(R)$, we need to analyze the necessary conditions under which the inequality $\alpha_{S_I}(R) \geq \alpha_{(I,T)}(R)$ holds for various fuzzy implications. Based on the theoretical analysis in Section 3.1, these conditions can be systematically derived. For conciseness, we present them in Table 3.
  \\{\bf Remark 8.} i.  According to Table 3, $R$ can be also referred to $\varepsilon$-transitive when $\alpha_{S_I}(R)\geq \varepsilon$.

 ii. As an indicator of $T$-transitivity, $\alpha_{S_I}(R)$ quantifies the degree of $T$-transitivity  from a logical view. While direct computation of  $\alpha_{S_I}(R)$ in Eq.(3) is non-trivial, Table 3 demonstrates that this degree can be effectively determined through $\alpha_{(I,T)}(R)$.\begin{center}
  \mbox{\bf{\small Table 3 The sufficient conditions for $\alpha_{S_I}(R)=\alpha_{(I,T)}(R)$}}
 \end{center}
 \begin{center}
\begin{tabular}{ccc}
    \toprule[1pt]
    Fuzzy implications &Sufficient conditions& Derivations\\
    \midrule
Continuous $I_{(S,N)}$& $T=(T_{\rm L})_\varphi$& Corollary 3.3\vspace{1mm}\\
  Continuous $I_T$ &$T$ is continuous &Corollary 3.8\vspace{1mm}\\
   $I_{(N,T_{\rm M},S)}(x,y)=\varphi^{-1}(\min(1,1-\varphi(x)+\varphi(y)))$& $T=(T_{\rm L})_\varphi$& Table 2\vspace{1mm}\\
   $I_g$ & $T=T_{\rm P}$ & Table 2\vspace{1mm}\vspace{1mm}\\
    $I_C$ & $T=T_{\rm P}$ & Table 2\vspace{1mm}\vspace{1mm}\\
    $\widetilde{I}_C$& $T=T_{\rm L}$ & Table 2\vspace{1mm}\vspace{1mm}\\
    \bottomrule[1pt]
  \end{tabular}
\end{center}\vspace{1.5mm}

 We therefore propose an efficient computational method for evaluating the $T$-transitivity degree of $R$, where $R=(r_{ij})_{n\times n}$ is  assumed to be a square matrix. This method consists of the following steps:

Step 1. Select a t-norm $T$ and a fuzzy implication $I$.

Step 2. Compute $a^k_{ij}=T(r_{ik},r_{kj})$ and compare $a^k_{ij}$ with $a^{k+1}_{ij}$ for all $i,j,l=1,2,\cdots,n$.

Step 3. Compute $\alpha_{ij}=I(a^n_{ij},r_{ij})$  for all $i,j,l=1,2,\cdots,n$.

Step 4. Select $\alpha=\max_{i,j}\alpha{ij}$ as the $T$-transitivity degree of $R$.

The detailed implementation process of the method is presented in Algorithm 1.
\\{\bf Algorithm 1} Calculating the degree of $T$-transitivity for fuzzy relation.

\textbf{Input:} The fuzzy relation $R=(r_{ij})_{n\times n}$.

\textbf{Output:} $\alpha$ ($T$-transitivity degree of $R$).
\\ 1. Select the t-norm $T$ and  fuzzy implication $I$.
\\ 2. \textbf{for} $(i=1, j=1, k=1;i++, j++, k++; i<n, j<n, k<n)$ \textbf{do}
\\3. \quad Compute $a^k_{ij}=T(r_{ik},r_{kj})$\vspace{1mm}
\\4. \quad \textbf{if} $a^k_{ij}>a^{k+1}_{ij}$ \textbf{then}\vspace{1mm}
\\5. \qquad $a^{k+1}_{ij}=a^{k}_{ij}$
\\ 6. \quad \textbf{end if}
\\7. \textbf{end for}
\\8. \textbf{for} each $i$ and $j$ \textbf{do}
\\9. \quad Compare $a^n_{ij}$ with $r_{ij}$
\\10. \quad \textbf{if} $a^n_{ij}>r_{ij}$ \textbf{then}
\\11. \qquad Compute $\alpha_{ij}=I(a^n_{ij},r_{ij})$
\\12. \quad \textbf{else}
\\13. \qquad $\alpha_{ij}=1$
\\14. \quad \textbf{end if}
\\15. \textbf{end for}
\\16. \textbf{for}  each $i$ \textbf{do}
\\17 \quad $j =1; j++; j<n$
\\18. \quad \textbf{if} $\alpha_{ij+1}<\alpha_{ij}$ \textbf{then}
\\19. \qquad $\alpha_{ij}=\alpha_{ij+1}$
\\20.  \quad \textbf{end if}
\\21. \textbf{end for}
\\22. \textbf{for}  each $j$ \textbf{do}
\\23. \quad $i=1; i++; i<n$
\\24. \quad \textbf{if}  $\alpha_{i+1j}<\alpha_{ij}$ \textbf{then}
\\25. \qquad $\alpha_{ij}=\alpha_{i+1j}$
\\26.  \quad \textbf{end if}
\\27. \textbf{end for}
\\28. $\alpha=\alpha_{ij}$.

The computational complexity and space complexity of Algorithm 1 are analyzed as follows.

i. Computational complexity. In step 2,  computing $a^k_{ij}=T(r_{ik},r_{kj})$ and comparing  $a^k_{ij}$ with $a^{k+1}_{ij}$ require  $2n^3$ operations. Similarly, step 3 also incurs $2n^3$ operations. Step 4 needs to compare $n$ times. Therefore, The algorithm's overall complexity is $\mathcal{O}(n^3)$, with the exact count being $4n^3+2n$ operations.

ii. Space complexity. In step 2,  storing $a_{ij}$ requires $n^2$ memory units.  Step 3 needs $n^2$ memory units to store $\alpha_{ij}$. Therefore, the space complexity of Algorithm 1  amounts to is $\mathcal{O}(n^2)$, with exact memory requirements of $2n^2+1$ units.
\\{\bf Example 3.12} Let $I(x,y)=\min\left(1,\frac{y}{x}\right)$ and $T=T_{\rm P}$. For a  fuzzy relation $R$ expressed as
$$R = \left(
\begin{array}{*{12}{c}}
 1.000 & 0.864 & 0.794 & 0.953 & 0.824 & 0.884 & 0.724 & 0.864 & 0.743 & 0.634 & 0.715 & 0.891 \\
0.984 & 1.000 & 0.662 & 0.954 & 0.883 & 0.793 & 0.713 & 0.662 & 0.862 & 0.565 & 0.743 & 0.752 \\
0.963 & 0.862 & 1.000 & 0.794 & 0.974 & 0.824 & 0.604 & 0.724 & 0.954 & 0.752 & 0.685 & 0.924 \\
0.537 & 0.593 & 0.717 & 1.000 & 0.824 & 0.974 & 0.824 & 0.743 & 0.884 & 0.576 & 0.629 & 0.653 \\
0.525 & 0.604 & 0.764 & 0.527 & 1.000 & 0.774 & 0.954 & 0.527 & 0.534 & 0.813 & 0.752 & 0.517 \\
0.462 & 0.882 & 0.593 & 0.536 & 0.704 & 1.000 & 0.624 & 0.624 & 0.824 & 0.928 & 0.413 & 0.814 \\
0.623 & 0.605 & 0.794 & 0.974 & 0.794 & 0.792 & 1.000 & 0.547 & 0.824 & 0.967 & 0.845 & 0.671 \\
0.927 & 0.817 & 0.727 & 0.864 & 0.703 & 0.862 & 0.704 & 1.000 & 0.713 & 0.752 & 0.576 & 0.908 \\
0.768 & 0.686 & 0.863 & 0.573 & 0.704 & 0.924 & 0.573 & 0.894 & 1.000 & 0.695 & 0.564 & 0.695 \\
0.723 & 0.634 & 0.785 & 0.845 & 0.962 & 0.745 & 0.693 & 0.576 & 0.767 & 1.000 & 0.482 & 0.427 \\
0.657 & 0.795 & 0.671 & 0.719 & 0.883 & 0.629 & 0.845 & 0.413 & 0.752 & 0.482 & 1.000 & 0.864 \\
0.891 & 0.752 & 0.924 & 0.653 & 0.517 & 0.614 & 0.671 & 0.908 & 0.695 & 0.876 & 0.864 & 1.000\\
\end{array}\right),$$
we determine that its degree of $T_{\rm P}$-transitivity  is $0.53$ using Algorithm 1. Moreover, $R$ is a 0.53-$T_P$-transitive fuzzy relation according to Definition 3.1.
  \section{Aggregation of $\varepsilon$-$T$-transitive fuzzy relations}
   \qquad  This section investigates the  aggregation of $\varepsilon$-$T$-transitive fuzzy relations. Let $R_i(i=1,\cdots, n)$ be a series of fuzzy relations defined on $U$. Given an $n$-ary aggregation function $G_n$, we define a fuzzy relation $R_{G_n}$ on $U$ as $R_{G_n}(x,y)=G_n(R_1(x,y),\cdots, R_n(x,y))$ for any $x,y\in U$\cite{Fodor}. A key question arises: how to select the aggregation function $G_n$ to guarantee $\varepsilon$-$T$-transitivity of $R_{G_n}$ when aggregating $\varepsilon$-$T$-transitive fuzzy relations $R_i(i=1,\cdots, n)$. To address this, we formally introduce the following definition.
   \\{\bf Definition 4.1} We say that the $n$-ary aggregation function $G_n$ aggregates a series of $\varepsilon$-$T$-transitive fuzzy relations $R_i(i=1,\cdots, n)$ if $R_{G_n}$ remains $\varepsilon$-$T$-transitive.

  Since the transitivity measure of $R$ is closely related to fuzzy implications, we will investigate the
preservation of $\varepsilon$-$T$-transitivity for a series of $\varepsilon$-$T$-transitive fuzzy relations based on the fuzzy implications discussed in Section 3.
   \\{\bf Proposition 4.2} Let  $I$ in $\alpha_{(I,T)}(R)$ be a continuous $(S,N)$-implication. $G_n$ aggregates a series of $\varepsilon$-$T$-transitive fuzzy relations $R_i(i=1,\cdots, n)$ $\Longleftrightarrow T(G_n(a_1,\cdots,a_n),G_n(b_1,\cdots,b_n))\leq \varphi^{-1}(\min(1,$ $\varphi(G_n(\varphi^{-1}(\max(0,\varphi(T(a_1,b_1))+\varphi(\varepsilon)-1)),\cdots,\varphi^{-1}(\max(0,$ $\varphi(T(a_n,b_n))+\varphi(\varepsilon)-1))))-\varphi(\varepsilon)+1))$ for any $a_i$ and $b_i(i=1,\cdots,n)\in [0,1]$.
   \\{\bf Proof.} $(\Longleftarrow)$ To guarantee that $R_{G_n}$ is $\varepsilon$-$T$-transitive, it needs to verify that the inequality $T(R_{G_n}(x,y),$ $R_{G_n}(y,z))\leq \varphi^{-1}(\min (1,1+\varphi(R_{G_n}(x,z))-\varphi(\varepsilon)))$ holds for all $x,y,z\in U$ according to Proposition 3.2. As $R_i(i=1,\cdots,n)$ is $\varepsilon$-$T$-transitive, we have $T(R_i(x,y),R_i(y,z))\leq \varphi^{-1}(\min (1,1+\varphi(R_i(x,z))-\varphi(\varepsilon)))$ for all $x,y,z\in U$ by Proposition 3.2. And then  $\varphi(T(R_i(x,y),$ $R_i(y,z)))\leq \varphi(\varphi^{-1}(\min (1,1+\varphi(R_i(x,z))-\varphi(\varepsilon))))=\min (1,1+\varphi(R_i(x,z))-\varphi(\varepsilon))$ holds for $i=1,\cdots,n$. This implies that $T(R_{G_n}(x,y),$ $R_{G_n}(y,z))\leq \varphi^{-1}(\min(1,\varphi(G_n(\varphi^{-1}(\max(0,\varphi(T(R_1(x,$ $y),R_1(y,z)))+\varphi(\varepsilon)-1)),\cdots,\varphi^{-1}(\max(0,$ $\varphi(T(R_1(x,y),R_1(y,z))))+\varphi(\varepsilon)-1))))-\varphi(\varepsilon)+1))\leq \varphi^{-1}(\min(1,\varphi(G_n(\varphi^{-1}(\max(0,\varphi(\varphi^{-1}(\min (1,$ $1+\varphi(R_1(x,z))-\varphi(\varepsilon)))+\varphi(\varepsilon)-1)),\cdots,\varphi^{-1}(\max(0,$ $\varphi(\varphi^{-1}(\min (1,1+\varphi(R_n(x,z))-\varphi(\varepsilon)))+\varphi(\varepsilon)-1))))-\varphi(\varepsilon)+1))\leq \varphi^{-1}(\min (1,1+\varphi(R_{G_n}(x,z))-\varphi(\varepsilon)))$. Thus,  $R_{G_n}$ is $\varepsilon$-$T$-transitive.

$(\Longrightarrow)$ For arbitrary $a_i$ and $b_i(i=1,\cdots,n)\in [0,1]$, let us define  a series of fuzzy relations $R_1, \cdots, R_n$ in the following: ${{R}_{i}}(x,x)={{R}_{i}}(y,y)={{R}_{i}}(z,z)=1$, $R_i(x,y)=R_i(y,x)=a_i$, $R_i(y,z)=R_i(z,y)=b_i$ and $R_i(x,z)=R_i(z,x)=\varphi^{-1}(\max(0,\varphi(T(a_i,b_i))+\varphi(\varepsilon)-1))(i=1,\cdots,n)$, where $x,y,z\in U$. Using Proposition 3.2, we can readily verify the $\varepsilon$-$T$-transitivity of each $R_i(i=1,\cdots,n)$. Since  $G_n$ aggregates  a series of $\varepsilon$-$T$-transitive fuzzy relations $R_i(i=1, \cdots, n)$, $R_{G_n}$ consequently maintains $\varepsilon$-$T$-transitivity. This means that $T(R_{G_n}(x,y),R_{G_n}(y,z))\leq \varphi^{-1}(\min (1,1+\varphi(R_{G_n}(x,z))-\varphi(\varepsilon)))$ holds for all $x,y,z\in U$ according to Proposition 3.2. Therefore, we can obtain $T(G_n(a_1,\cdots,a_n),$ $G_n(b_1,\cdots,b_n))=T(G_n(R_1(x,y),\cdots,R_n(x,y)),$ $G_n(R_1(y,z),\cdots,R_n(y,z)))\leq \varphi^{-1}(\min (1,1+\varphi(R_{G_n}$ $(x,z))-\varphi(\varepsilon)))=\varphi^{-1}(\min(1,\varphi(G_n(\varphi^{-1}$ $(\max(0,\varphi(T(a_1,b_1))+\varphi(\varepsilon)-1)),\cdots,\varphi^{-1}(\max(0,\varphi(T(a_n, b_n))+\varphi(\varepsilon)-1))))-\varphi(\varepsilon)+1))$.$\hfill\square$
\\{\bf Proposition 4.3} Let  $I$ in $\alpha_{(I,T)}(R)$ be a R-implication generated by a left-continuous $T'$.  $G_n$ aggregates a series of $\varepsilon$-$T$-transitive fuzzy relations $R_i(i=1,  \cdots, n) \Longleftrightarrow T'(\varepsilon, T(G_n(a_1,\cdots,$ $a_n),G_n(b_1,\cdots,b_n)))\leq G_n(T'(\varepsilon,T(a_1,b_1)), \cdots,T'(\varepsilon,T(a_n,b_n)))$ holds for all $a_i$ and $b_i(i=1,\cdots,n)\in [0,1]$.
   \\{\bf Proof.} $(\Longleftarrow)$ To ensure that $R_{G_n}$ is $\varepsilon$-$T$-transitive, it needs to verify that the inequality $T'(\varepsilon,T(R_{G_n}(x,y),$ $R_{G_n}(y,z)))\leq R_{G_n}(x,z)$ holds for all $x,y,z\in U$ by Proposition 3.6. Since $R_i(i=1,\cdots,n)$ is $\varepsilon$-$T$-transitive, $T'(\varepsilon,T(R_i(x,y),R_i(y,z)))\leq R_i(x,z)$ holds for all $x,y,z\in U$ according to Proposition 3.6. Thus, we have  $T'(\varepsilon,T(R_{G_n}(x,y),R_{G_n}(y,z)))=T'(\varepsilon,T(G_n(R_1(x,y),\cdots,R_n(x,y)),G_n(R_1(y,z),\cdots,R_n(y,z)))\leq G_n(T'(\varepsilon,T(R_1(x,y),R_1(y,z))),$ $ \cdots,T'(\varepsilon,T(R_n(x,y),R_n(y,z))))\leq {G_n}(R_1(x,z),\cdots,R_n(x,z))= R_{G_n}(x,z)$.

$(\Longrightarrow)$ For arbitrary $a_i$ and $b_i(i=1,\cdots,n)\in [0,1]$, let us define  a series of fuzzy relations $R_1, \ldots, R_n$ in the following: ${{R}_{i}}(x,x)={{R}_{i}}(y,y)={{R}_{i}}(z,z)=1$, $R_i(x,y)=R_i(y,x)=a_i$, $R_i(y,z)=R_i(z,y)=b_i$ and $R_i(x,z)=R_i(z,x)=T'(T(a_i,b_i),\varepsilon)(i=1,\cdots,n)$, where $x,y,z\in U$. According to  Proposition 3.6,  $R_1,\cdots,R_n$ is a series of $\varepsilon$-$T$-transitive fuzzy relations on $U$. Since  $G_n$ aggregates $R_1, \cdots, R_n$, we have the fact that $R_{G_n}$ is  $\varepsilon$-$T$-transitive, too. Therefore,  $T'(\varepsilon, T(R_{G_n}(x,y),R_{G_n}(y,z)))\leq R_{G_n}(x,z)$ holds by Proposition 3.6. This implies that
$T'(\varepsilon, T(G_n(a_1,\cdots,a_n),G_n(b_1,\cdots,b_n)) =T'(\varepsilon,T({G_n}(R_1(x,y),\cdots,R_n(x,y)),{G_n}(R_1(y,z),\cdots,$ $R_n(y,z))))\leq {G_n}(R_1(x,z),\cdots,R_n(x,z))={G_n}(T'(\varepsilon,T(a_1,b_1)), \cdots,$ $T'(\varepsilon,T(a_n,b_n)))$ holds. $\hfill\square$
\\{\bf Proposition 4.4} Let $I$ in $\alpha_{(I,T)}(R)$ be a $T$-power implication with an additive generator $t$. $G_n$ aggregates a series of $\varepsilon$-$T$-transitive fuzzy relations $R_1, \cdots, R_n \Longleftrightarrow T(G_n(a_1,\cdots,a_n),G_n(b_1,$\vspace{1mm} $\cdots,b_n))\leq t^{-1}\left(\varepsilon t\left(G_n\left(t^{-1}\left(\frac{t(T(a_1,b_1))}{\varepsilon}\right),\cdots, t^{-1}\left(\frac{t(T(a_n,b_n))}{\varepsilon}\right)\right)\right)\right)$ for any $a_i$ and $b_i(i=1,$\vspace{1mm} $\cdots,n)\in [0,1]$.
   \\{\bf Proof.} $(\Longleftarrow)$ To verify that $R_{G_n}$ is $\varepsilon$-$T$-transitive, we need to ensure that the inequality $T(R_{G_n}(x,y),$ $R_{G_n}(y,z))\leq t^{-1}(\varepsilon t(R_{G_n}(x,z)))$ holds for all $x, y, z\in U$ according to Proposition 3.12. Since $R_i(i=1,\cdots,n)$ is $\varepsilon$-$T$-transitive, Proposition 3.12 implies that $T(R_i(x,y),R_i(y,z))$ $\leq t^{-1}(\varepsilon t(R_i(x,z)))$ holds for all $x, y, z\in U$. Therefore, we can obtain $T(R_{G_n}(x,y),R_{G_n}(y,z))\leq$\vspace{1mm} $ t^{-1}\left(\varepsilon t\left( G_n\left(t^{-1}\left(\frac{t(T(R_1(x,y),R_1(y,z)))}{\varepsilon}\right)\right.\right.\right.,$ $\left.\left.\left.\cdots, t^{-1}\left(\frac{t(T(R_n(x,y),R_n(y,z)))}{\varepsilon}\right)\right)\right)\right)$\vspace{1mm} $\leq t^{-1}\left(\varepsilon t\left(G_n\left(t^{-1}\right.\right.\right.$\vspace{1mm} $\left(\frac{t(t^{-1}(\varepsilon t(R_1(x,z))))}{\varepsilon}\right),\cdots,\left.\left.\left. t^{-1}\left(\frac{t(t^{-1}(\varepsilon t(R_n(x,z))))}{\varepsilon}\right)\right)\right)\right)=$ $t^{-1}(\varepsilon t(R_{G_n}(x,z)))$.\vspace{1mm}

$(\Longrightarrow)$ For any given $a_i$ and $b_i(i=1,\cdots,n)\in [0,1]$, let us define  a series of fuzzy relations $R_1, \cdots, R_n$ in the following: ${{R}_{i}}(x,x)={{R}_{i}}(y,y)={{R}_{i}}(z,z)=1$, $R_i(x,y)=R_i(y,x)=a_i$,\vspace{1mm} $R_i(y,z)=R_i(z,y)=b_i$ and $R_i(x,z)=R_i(z,x)=t^{-1}\left(\frac{t(T(a_i,b_i))}{\varepsilon}\right)$. It is not difficult to\vspace{1mm} verify that $R_i\ (i=1,\ldots,n)$ is a series of $\varepsilon$-$T$-transitive fuzzy relations by Proposition 3.12. As  $G_n$ aggregates  a series of $\varepsilon$-$T$-transitive fuzzy relations $R_1, R_2, \ldots, R_n$, $R_{G_n}$ is  $\varepsilon$-$T$-transitive, too. This means that $T(R_{G_n}(x,y),R_{G_n}(y,z))\leq t^{-1}(\varepsilon t(R_{G_n}(x,z)))$ holds for all $x,y,z\in U$ according to Proposition 3.12. Therefore, we have $T(G_n(a_1,\cdots,a_n),G_n(b_1$ $,\cdots,b_n))=T(G_n(R_1(x,y),\cdots,R_n(x,y)),G_n$ $(R_1$ $(y,z),\cdots,R_n(y,z)))\leq t^{-1}(\varepsilon t(R_{G_n}(x,z)))=t^{-1}\left(\varepsilon t\left(G_n\right.\right.$\vspace{1mm} $\left(t^{-1}\left(\frac{t(T(a_1,b_1))}{\varepsilon}\right),\cdots, \left.\left. t^{-1}\left(\frac{t(T(a_n,b_n))}{\varepsilon}\right)\right)\right)\right)$.$\hfill\square$\vspace{1mm}

For others fuzzy implications, analogous necessary and sufficient conditions can be established to characterize when $G_n$ preserves $\varepsilon$-$T$-transitivity. These conditions are systematically summarized in the following table.
 \begin{center}
  \mbox{\bf{\small Table 4 The necessary and sufficient conditions for $\varepsilon$-$T$-transitivity of $R_{G_n}$ }}\vspace{-4mm}
 \end{center}
 \begin{center}
\begin{tabular}{cr}
    \toprule[1pt]
    Fuzzy implications&  Necessary and sufficient conditions\\
    \midrule
  $I_{(N,T_{\rm M},S)}$ & $T(G_n(a_1,\cdots,a_n),G_n(b_1,\cdots,b_n))\leq N^{-1}(\varphi^{-1}(\max(0,\varphi(\varepsilon)-\varphi(G_n(\varphi^{-1}(\max(0,$\\
  &$\varphi(\varepsilon)-\varphi(N(T(a_1,b_1))))),\cdots,\varphi^{-1}(\max(0,\varphi(\varepsilon)-\varphi(N(T(a_n,b_n))))))))))$\vspace{1mm}\\
  $I_{(N,T_{\rm P},S)}$ & $T(G_n(a_1,\cdots,a_n),G_n(b_1,\cdots,b_n))\leq\min\left(1,\right.$\\
  &$\left.\frac{G_n(M(a_1,b_1),\cdots, M(a_n, b_n))+\sqrt{G_n^2(M(a_1,b_1),\cdots, M(a_n,b_n))+4(1-\varepsilon)}}{2}\right)$\vspace{1mm}\\
   $I_g$ & $T(G_n(a_1,\cdots,a_n),G_n(b_1,\cdots,b_n))\leq G_n(\varepsilon T(a_1,b_1),\cdots,$ $\varepsilon T(a_n,$ $b_n))$\vspace{1mm}\\
   $I_C$ & $T(G_n(a_1,\cdots,a_n),G_n(b_1,\cdots,b_n))\leq G_n(\varepsilon T(a_1,b_1),$ $\cdots,\varepsilon T(a_n,b_n))$\vspace{1mm}\\
    $\widetilde{I}_C$& $T(G_n(a_1,\cdots,a_n),G_n(b_1,\cdots,b_n))\leq \min(1,1-\varepsilon+G_n(\max(0,$\\
    &$T(a_1,b_1)+\varepsilon-1),\cdots, \max(0,T(a_1,b_1)+\varepsilon-1)))$ \vspace{1mm}\\
    \bottomrule[1pt]
  \end{tabular}
\end{center}
Where $M(a_i, b_i)=\max\left(0,T(a_i,b_i)+\frac{\varepsilon-1}{T(a_i,b_i)}\right)$.\vspace{1mm}

In the final part of this section, we examine the preservation of $\varepsilon$-$T$-transitivity for the aggregation functions presented in Example 2.2.
\\{\bf Example 4.5} We examine two fuzzy relations constructed by the Max and Min aggregation functions:
$$R_{\rm Max}(x,y)=\max_{i=1}^nR_i(x,y)\quad {\rm and} \quad R_{\rm Min}(x,y)=\min_{i=1}^nR_i(x,y).$$

It is straightforward  to verify that the Max and Min aggregation functions  meet the conditions specified in Propositions 4.2-4.4 and Table 4. Consequently, these functions preserve  $\varepsilon$-$T$-transitivity, implying that $R_{\rm Max}$ and $R_{\rm Min}$ are also $\varepsilon$-$T$-transitive.
\\{\bf Example 4.6} Suppose that $I_g$ is employed to measure the transitivity degree of a fuzzy relation.  Let $R_i(i=1,\cdots,n)$ be  a series of  $\varepsilon$-$T_{\rm P}$-transitive fuzzy relations on $U$. With the WQAM\vspace{1mm} $G_{f,\omega}(x_1,\cdots,x_n)=\exp\left(\mathop{\sum}\limits_{i=1}^n\omega_i\log x_i\right)$, we define a fuzzy relation $R_{G_{f,\omega}}$ as
$$R_{G_{f,\omega}}(x,y)=\exp\left(\sum_{i=1}^n\omega_i\log R_i(x,y)\right).$$

It can be verified that the WQAM  $G_{f,\omega}$ fulfills the conditions specified in Table 4. This implies that the WQAM  $G_{f,\omega}$ preserves  $\varepsilon$-$T$-transitivity, that is, $R_{G_{f,\omega}}$ is also $\varepsilon$-$T$-transitive.
   \section{Some applications of $\varepsilon$-$T$-transitive fuzzy relation}
   \qquad This section will consider the applications of $\varepsilon$-$T$-transitive fuzzy relation in approximate reasoning and decision-making. Especially,  the equivalence fuzzy relation will be replaced by a $\varepsilon$-equivalence fuzzy relation to make a decision.
\subsection{Approximate reasoning with $\varepsilon$-trnasitive fuzzy relation}
 \qquad  Let us consider the well-known generalization of hypothetical syllogism (GHS) intuitively
expressed as
\begin{verse}
	\hspace{37mm}	Premise 1: IF $x$ is $A$ THEN $y$ is $B$\\
	\hspace{37mm}  Premise 2: IF $y$ is $B'$ THEN $z$ is $C$\\
	\hspace{35mm}   \rule[2pt]{6cm}{0.05em}\\
	\hspace{37mm}   Conclusion: IF $x$ is $A$ THEN $z$ is $C',$
\end{verse}
where $A, B$ and $B'$, $C$ and $C'$ are five fuzzy sets on the universes $U, V$ and $W$, respectively. Suppose that $R(x, y)$ and $R(y, z)$ are used to interpret the truth-values of Promises 1 and 2, respectively. We further choose a t-norm
$T$ and a fuzzy implication $I$ to respectively interpret the
fuzzy logical connectives ``conjunction" and ``conditional statement". Then, the
truth-value\vspace{1mm} of conclusion can be expressed as $I\left(\mathop{\sup}\limits_{y\in V}T(R(x, y),R(y',z)), R(x,z')\right)$. In this expression, we\vspace{1mm} expect that the closer $R(x, z')$ is to $R(x,z)$, the closer $y$ is to $y'$. Especially,\vspace{1mm} $R(x, z')=R(x, z)$ when $y=y'$. This means that  $I\left(\mathop{\sup}\limits_{y\in V}T(R(x, y),R(y',z)), R(x,z')\right)\geq \varepsilon$\vspace{1mm} is expected to hold if $R(x, y) \geq \varepsilon$
and $R(y, x)\geq \varepsilon$. Obviously, it is difficult to achieve this expectation for a non-transitive fuzzy
relation $R$. However, we can choose some $\varepsilon$-$T$-transitive fuzzy relations to express the truth-values of
Promises 1 and 2. Then, a suitable t-norm can be selected according to the conditions of Propositions 3.2, 3.6, 3.9  and Table 1 in order to achieve this expectation. For example, let
the truth-values of Promises 1 and 2 are at least 0.8. If the R-implication
$I(x, y) =\min(1, \frac{y}{x})$ is used to interpreted the logical connective ``conditional statement", we\vspace{1mm}
can ensure that the truth-value of conclusion  $I\left(\mathop{\sup}\limits_{y\in V}T(R(x, y),R(y',z)), R(x,z')\right)\geq 0.8$ holds\vspace{1mm} when $T= T_{\rm P}$.

Further, let $T$ be a continuous Archimedean t-norm, that is, $T(x, y)=t^{-)}(\min(t(0),t(x)+t(y)))$ with an additive generator $t$. For an $\varepsilon$-$T$-transitive relation $R$, we define a pseudo-metric on $U$ as $d(x, y) = t(R(x, y))$. Indeed, the inequality  $d(x, z) = t(R(x, z))\leq t(\varepsilon) + t(R(x, y)) + t(R(y, z)) =
t(\varepsilon) + d(x, y) + d(y, z)$ holds for any $x, y, z\in U$. Obviously, $d$ is a metric when $R$ becomes an equivalent fuzzy relation. Therefore, the pseudo-metric $d$ can be used to measure the distance to a certain extent. Let us return to the discussion of GHS. Suppose that $d(x, y)\leq \delta$ and $d(y, z)\leq \delta$. Then, we have $d(x, z)\leq t(\varepsilon)+2\delta$. This means that the closer $R(x, z')$) is to $R(x, z)$ when $y$ is closer to $y'$.
  \subsection{Clustering with $\varepsilon$-equivalent fuzzy relation}
\qquad This subsection extends the concept of equivalence fuzzy relation to $(\varepsilon_1,\varepsilon_2)$-equivalence fuzzy relation to meet practical requirements.  According to Definition 3.1 in \cite{Dan},  $R$ is referred as a $\varepsilon$-symmetric fuzzy relation if its degree of symmetry $\alpha_s(R)\geq \varepsilon$. Building on this concept, we define a $(\varepsilon_1,\varepsilon_2)$-equivalence fuzzy relation as follows.
\\{\bf Definition 5.1} Let $R$ be a reflexive fuzzy relation on $U$ and $\varepsilon_1,\varepsilon_2 \in[0,1]$. We say that $R$ is $(\varepsilon_1,\varepsilon_2)$-equivalent when it is both $\varepsilon_1$-symmetric and $\varepsilon_2$-$T$-transitive. Especially, $R$ is $\varepsilon$-equivalent if it is reflexive, symmetric and $\varepsilon$-$T$-transitive.

   Similar to the application of equivalence fuzzy relation, we will cluster 18 objects with a $\varepsilon$-equivalence fuzzy relation as shown in the following example.
\\{\bf Example 5.2}\cite{Li} As one of the important indicators, the vibration state can be used to assess the working condition of turbine-generators. To analyze vibration faults in turbine-generators, oil whip, unbalance, and misalignment are selected as three typical vibration faults. Sample data are collected based on the frequency domain characteristics of vibration signals from the following five frequency bands:  $<0.4f$, 0.4-0.5$f$, $1f$, $2f$, and $\geq 3f$, where $f$ denotes the rotational frequency. The normalized amplitude components can be presented in Table 5.
$$\mbox{\bf{\small Table\ 5\ The fault sample of a 50 MW turbine-generator}}$$
\begin{center}
  \begin{tabular}{ccccccl}
    \toprule
    No. & \multicolumn{5}{c}{Sample input} & Fault type \\
    \cmidrule(lr){2-6}
    & $<$0.4$f$ & 0.4-0.5$f$ & 1$f$ & 2$f$ & $\geq$3$f$ & \\
    \midrule
    1  & 0.052  & 0.783  & 0.225 & 0.036 & 0.013 & Oil whirl    \\
    2  & 0.232  & 0.975  & 0.314 & 0.056 & 0.014 & Oil whirl    \\
    3  & 0.161  & 0.925  & 0.285 & 0.023 & 0.016 & Oil whirl    \\
    4  & 0.106  & 0.858  & 0.230 & 0.017 & 0.028 & Oil whirl    \\
    5  & 0.079  & 0.819  & 0.201 & 0.016 & 0.012 & Oil whirl    \\
    6  & 0.028  & 0.061  & 0.980 & 0.225 & 0.057 & Unbalance    \\
    7  & 0.045  & 0.022  & 1.000 & 0.316 & 0.065 & Unbalance    \\
    8  & 0.010  & 0.054  & 0.875 & 0.183 & 0.073 & Unbalance    \\
    9  & 0.015  & 0.032  & 0.923 & 0.219 & 0.037 & Unbalance    \\
    10 & 0.023  & 0.025  & 0.758 & 0.115 & 0.019 & Unbalance    \\
    11 & 0.033  & 0.037  & 0.386 & 0.531 & 0.230 & Misalignment \\
    12 & 0.017  & 0.023  & 0.397 & 0.458 & 0.103 & Misalignment \\
    13 & 0.012  & 0.039  & 0.427 & 0.496 & 0.175 & Misalignment \\
    14 & 0.021  & 0.017  & 0.298 & 0.403 & 0.132 & Misalignment \\
    15 & 0.017  & 0.056  & 0.483 & 0.599 & 0.301 & Misalignment \\
    \bottomrule
  \end{tabular}
\end{center}\vspace{2mm}

Table 6 lists three datasets capturing the frequency-domain characteristics of vibration signals from analyzed samples.
$$\mbox{\bf{\small Table\ 6\ Test data of  a 50 MW turbine-generator}}$$
\begin{center}
  \begin{tabular}{lccccc}
    \toprule
    No. & $<$0.4$f$ & 0.4$\sim$0.5$f$ & 1$f$ & 2$f$ & $\geq$3$f$ \\
    \midrule
    16  & 0.161 & 0.753          & 0.128 & 0.006 & 0.003  \\
    17  & 0.017 & 0.053          & 0.750 & 0.252 & 0.107  \\
    18  & 0.026 & 0.043          & 0.357 & 0.517 & 0.098  \\
    \bottomrule
  \end{tabular}
\end{center}\vspace{2mm}

To diagnose the vibration fault of the turbine-generator using the three datasets from the checked samples in Table 5, the fault sample data is combined with the checked sample data. The fuzzy relation $R$ can then be obtained using the cosine of included angle method as follows \cite{Hu}.
\[R = \resizebox{1\textwidth}{!}{
$\left(
\begin{array}{*{18}{c}}
1.0000 & 0.9934 & 0.9973 & 0.9989 & 0.9992 & 0.6688 & 0.6493 & 0.6683 & 0.6558 & 0.6560 & 0.6225 & 0.6253 & 0.6328 & 0.6151 & 0.6321 & 0.9914 & 0.6694 & 0.6297\vspace{1mm} \\
0.9934 & 1.0000 & 0.9989 & 0.9963 & 0.9943 & 0.6823 & 0.6647 & 0.6804 & 0.6688 & 0.6699 & 0.6349 & 0.6376 & 0.6437 & 0.6277 & 0.6423 & 0.9945 & 0.6822 & 0.6419\vspace{1mm}  \\
0.9973 & 0.9989 & 1.0000 & 0.9991 & 0.9979 & 0.6748 & 0.6554 & 0.6736 & 0.6612 & 0.6631 & 0.6208 & 0.6236 & 0.6305 & 0.6133 & 0.6296 & 0.9953 & 0.6738 & 0.6271\vspace{1mm}  \\
0.9989 & 0.9963 & 0.9991 & 1.0000 & 0.9997 & 0.6583 & 0.6385 & 0.6578 & 0.6447 & 0.6460 & 0.6117 & 0.6121 & 0.6204 & 0.6032 & 0.6211 & 0.9957 & 0.6585 & 0.6167 \vspace{1mm} \\
0.9992 & 0.9943 & 0.9979 & 0.9997 & 1.0000 & 0.6482 & 0.6281 & 0.6476 & 0.6347 & 0.6359 & 0.6030 & 0.6041 & 0.6120 & 0.5950 & 0.6122 & 0.9954 & 0.6481 & 0.6096 \vspace{1mm} \\
0.6688 & 0.6823 & 0.6748 & 0.6583 & 0.6482 & 1.0000 & 0.9979 & 0.9997 & 0.9997 & 0.9981 & 0.8652 & 0.9023 & 0.8961 & 0.8740 & 0.8763 & 0.6126 & 0.9961 & 0.8690 \vspace{1mm} \\
0.6493 & 0.6647 & 0.6554 & 0.6385 & 0.6281 & 0.9979 & 1.0000 & 0.9968 & 0.9983 & 0.9936 & 0.8892 & 0.9245 & 0.9177 & 0.8986 & 0.8981 & 0.5937 & 0.9980 & 0.8941 \vspace{1mm} \\
0.6683 & 0.6804 & 0.6736 & 0.6578 & 0.6476 & 0.9997 & 0.9968 & 1.0000 & 0.9992 & 0.9982 & 0.8623 & 0.8979 & 0.8930 & 0.8702 & 0.8744 & 0.6109 & 0.9958 & 0.8635 \vspace{1mm} \\
0.6558 & 0.6688 & 0.6612 & 0.6447 & 0.6347 & 0.9997 & 0.9983 & 0.9992 & 1.0000 & 0.9982 & 0.8643 & 0.9032 & 0.8960 & 0.8743 & 0.8750 & 0.5985 & 0.9954 & 0.8698 \vspace{1mm} \\
0.6560 & 0.6699 & 0.6631 & 0.6460 & 0.6359 & 0.9981 & 0.9936 & 0.9982 & 0.9982 & 1.0000 & 0.8360 & 0.8773 & 0.8697 & 0.8459 & 0.8480 & 0.6004 & 0.9893 & 0.8406 \vspace{1mm} \\
0.6225 & 0.6349 & 0.6208 & 0.6117 & 0.6030 & 0.8652 & 0.8892 & 0.8623 & 0.8643 & 0.8360 & 1.0000 & 0.9911 & 0.9968 & 0.9983 & 0.9989 & 0.5794 & 0.9035 & 0.9916 \vspace{1mm} \\
0.6253 & 0.6376 & 0.6236 & 0.6121 & 0.6041 & 0.9023 & 0.9245 & 0.8979 & 0.9032 & 0.8773 & 0.9911 & 1.0000 & 0.9977 & 0.9965 & 0.9891 & 0.5769 & 0.9321 & 0.9968 \vspace{1mm} \\
0.6328 & 0.6437 & 0.6305 & 0.6204 & 0.6120 & 0.8961 & 0.9177 & 0.8930 & 0.8960 & 0.8697 & 0.9968 & 0.9977 & 1.0000 & 0.9984 & 0.9966 & 0.5844 & 0.9292 & 0.9943\vspace{1mm}  \\
0.6151 & 0.6277 & 0.6133 & 0.6032 & 0.5950 & 0.8740 & 0.8986 & 0.8702 & 0.8743 & 0.8459 & 0.9983 & 0.9965 & 0.9984 & 1.0000 & 0.9959 & 0.5705 & 0.9098 & 0.9968 \vspace{1mm} \\
0.6321 & 0.6423 & 0.6296 & 0.6211 & 0.6122 & 0.8763 & 0.8981 & 0.8744 & 0.8750 & 0.8480 & 0.9989 & 0.9891 & 0.9966 & 0.9959 & 1.0000 & 0.5860 & 0.9138 & 0.9869 \vspace{1mm} \\
0.9914 & 0.9945 & 0.9953 & 0.9957 & 0.9954 & 0.6126 & 0.5937 & 0.6109 & 0.5985 & 0.6004 & 0.5794 & 0.5769 & 0.5844 & 0.5705 & 0.5860 & 1.0000 & 0.6124 & 0.5860 \vspace{1mm} \\
0.6694 & 0.6822 & 0.6738 & 0.6585 & 0.6481 & 0.9961 & 0.9980 & 0.9958 & 0.9954 & 0.9893 & 0.9035 & 0.9321 & 0.9292 & 0.9098 & 0.9138 & 0.6124 & 1.0000 & 0.9028 \vspace{1mm} \\
0.6297 & 0.6419 & 0.6271 & 0.6167 & 0.6096 & 0.8690 & 0.8941 & 0.8635 & 0.8698 & 0.8406 & 0.9916 & 0.9968 & 0.9943 & 0.9968 & 0.9869 & 0.5860 & 0.9028 & 1.0000 \vspace{1mm} \\
\end{array}
\right)$}\]

Let $T=T_{\rm P}$ and $I(x,y)=\min\left(1,\frac{y}{x}\right)$. According Algorithm 1,  The degree of transitivity of this fuzzy relation is $\alpha_{(I,T)}(R)= 0.8981$. By Definition 5.1, $R$ is  a 0.8981-equivalence fuzzy relation. Consequently, $R_\lambda$ can be regarded as an equivalence relation for some $\lambda\in [0,1]$\cite{Hu}. Thus, the data of fault samples and
checked samples can be clustered with $R_\lambda$ as shown in Table 7.
$$\mbox{\bf{\small Table\ 7\ Clustering with $R$ at different $\lambda$}}$$
  \begin{center}
  \begin{tabular}{cc}
    \toprule[1pt]
    $\lambda$ &  Clustering of the fault samples and check samples \\
    \midrule
    0.91 & $\{1, 2, 3, 4, 5, 16\}, \{6,7, 8, 9, 10,12,13,15,17\}, \{7,11, 12, 13, 14, 15, 17,18\}$\vspace{1mm}\\
   0.92 & $\{1, 2, 3, 4, 5, 16\}, \{6,7, 8, 9, 10,12,13,17\}, \{7,11, 12, 13, 14, 15, 17,18\}$\vspace{1mm}\\
   0.93 & $\{1, 2, 3, 4, 5, 16\}, \{6,7, 8, 9, 10,12,17\}, \{11, 12, 13, 14, 15, 17,18\}$\vspace{1mm}\\
    0.94 & $\{1, 2, 3, 4, 5, 16\}, \{6,7, 8, 9, 10,17\}, \{11, 12, 13, 14, 15,18\}$\vspace{1mm}\\
    0.95 & $\{1, 2, 3, 4, 5, 16\}, \{6,7, 8, 9, 10,17\}, \{11, 12, 13, 14, 15,18\}$\vspace{1mm}\\
    0.96 & $\{1, 2, 3, 4, 5, 16\}, \{6,7, 8, 9, 10,17\}, \{11, 12, 13, 14, 15, 18\}$\vspace{1mm}\\
    0.97 & $\{1, 2, 3, 4, 5, 16\}, \{6,7, 8, 9, 10,17\}, \{11, 12, 13, 14, 15,18\}$\vspace{1mm}\\
   0.98 & $\{1, 2, 3, 4, 5, 16\}, \{6,7, 8, 9, 10,17\}, \{11, 12, 13, 14, 15, 18\}$\vspace{1mm}\\
    0.99 & $\{1, 2, 3, 4, 5, 16\}, \{6,7, 8, 9, 10,17\}, \{11, 12, 13, 14, 15, 18\}$\vspace{1mm}\\
    1.00&\{1\},\{2\},\{3\},\{4\},\{5\},\{6\},\{7\},\{8\},\{9\},\{10\},\{11\},\{12\},\{13\},\{14\},\{15\},\{16\},\{17\},\{18\}\vspace{1mm}\\
    \bottomrule[1pt]
  \end{tabular}
\end{center}

Obviously, both fault samples and checked samples can be clustered into three distinct categories when  $\lambda$ varies from 0.9098 to 0.9914. This suggests that the checked samples 16, 17 and 18 can be diagnosed as oil whip, unbalance and  misalignment, respectively. These diagnostic results align perfectly with the actual conditions reported in \cite{Li}.

Let us reconsider Example 5.2 with the classical fuzzy clustering method. It can be observed that $R$ is not $T_{\rm P}$-transitive. Therefore, we need to compute the $T_{\rm P}$-transitive closure of $R$ using  Algorithm 1 from Ref.\cite{Naessens}. The $T_{\rm P}$-transitive closure $\overline{R}^T$ can be obtained as follows.\vspace{1mm}
\[
\overline{R}^T=\resizebox{1\textwidth}{!}{
$\left(
\begin{array}{*{18}{c}}
1.0000 & 0.9969 & 0.9980 & 0.9989 & 0.9992 & 0.6802 & 0.6788 & 0.6800 & 0.6800 & 0.6789 & 0.6397 & 0.6402 & 0.6417 & 0.6407 & 0.6403 & 0.9946 & 0.6801 & 0.6399\vspace{1mm}\\
0.9969 & 1.0000 & 0.9989 & 0.9980 & 0.9977 & 0.6823 & 0.6809 & 0.6821 & 0.6821 & 0.6810 & 0.6416 & 0.6422 & 0.6437 & 0.6427 & 0.6423 & 0.9945 & 0.6822 & 0.6419 \vspace{1mm}\\
0.9980 & 0.9989 & 1.0000 & 0.9991 & 0.9988 & 0.6815 & 0.6802 & 0.6813 & 0.6813 & 0.6803 & 0.6409 & 0.6415 & 0.6430 & 0.6420 & 0.6416 & 0.9953 & 0.6814 & 0.6412 \vspace{1mm}\\
0.9989 & 0.9980 & 0.9991 & 1.0000 & 0.9997 & 0.6809 & 0.6796 & 0.6807 & 0.6807 & 0.6796 & 0.6404 & 0.6409 & 0.6424 & 0.6414 & 0.6410 & 0.9957 & 0.6808 & 0.6406 \vspace{1mm}\\
0.9992 & 0.9977 & 0.9988 & 0.9997 & 1.0000 & 0.6807 & 0.6794 & 0.6805 & 0.6805 & 0.6794 & 0.6402 & 0.6407 & 0.6422 & 0.6412 & 0.6408 & 0.9954 & 0.6806 & 0.6404 \vspace{1mm}\\
0.6802 & 0.6823 & 0.6815 & 0.6809 & 0.6807 & 1.0000 & 0.9980 & 0.9997 & 0.9997 & 0.9981 & 0.9236 & 0.9285 & 0.9263 & 0.9252 & 0.9232 & 0.6785 & 0.9961 & 0.9255 \vspace{1mm}\\
0.6788 & 0.6809 & 0.6802 & 0.6796 & 0.6794 & 0.9980 & 1.0000 & 0.9977 & 0.9983 & 0.9965 & 0.9254 & 0.9302 & 0.9281 & 0.9270 & 0.9249 & 0.6772 & 0.9980 & 0.9273 \vspace{1mm}\\
0.6800 & 0.6821 & 0.6813 & 0.6807 & 0.6805 & 0.9997 & 0.9977 & 1.0000 & 0.9994 & 0.9982 & 0.9234 & 0.9282 & 0.9261 & 0.9249 & 0.9229 & 0.6783 & 0.9958 & 0.9252 \vspace{1mm}\\
0.6800 & 0.6821 & 0.6813 & 0.6807 & 0.6805 & 0.9997 & 0.9983 & 0.9994 & 1.0000 & 0.9982 & 0.9238 & 0.9287 & 0.9265 & 0.9254 & 0.9234 & 0.6783 & 0.9963 & 0.9257 \vspace{1mm}\\
0.6789 & 0.6810 & 0.6803 & 0.6796 & 0.6794 & 0.9981 & 0.9965 & 0.9982 & 0.9982 & 1.0000 & 0.9222 & 0.9270 & 0.9249 & 0.9237 & 0.9217 & 0.6773 & 0.9945 & 0.9240 \vspace{1mm}\\
0.6397 & 0.6416 & 0.6409 & 0.6404 & 0.6402 & 0.9236 & 0.9254 & 0.9234 & 0.9238 & 0.9222 & 1.0000 & 0.9948 & 0.9968 & 0.9983 & 0.9989 & 0.6381 & 0.9273 & 0.9951 \vspace{1mm}\\
0.6402 & 0.6422 & 0.6415 & 0.6409 & 0.6407 & 0.9285 & 0.9302 & 0.9282 & 0.9287 & 0.9270 & 0.9948 & 1.0000 & 0.9977 & 0.9965 & 0.9943 & 0.6387 & 0.9321 & 0.9968 \vspace{1mm}\\
0.6417 & 0.6437 & 0.6430 & 0.6424 & 0.6422 & 0.9263 & 0.9281 & 0.9261 & 0.9265 & 0.9249 & 0.9968 & 0.9977 & 1.0000 & 0.9984 & 0.9966 & 0.6402 & 0.9300 & 0.9952 \vspace{1mm}\\
0.6407 & 0.6427 & 0.6420 & 0.6414 & 0.6412 & 0.9252 & 0.9270 & 0.9249 & 0.9254 & 0.9237 & 0.9983 & 0.9965 & 0.9984 & 1.0000 & 0.9972 & 0.6391 & 0.9288 & 0.9968 \vspace{1mm}\\
0.6403 & 0.6423 & 0.6416 & 0.6410 & 0.6408 & 0.9232 & 0.9249 & 0.9229 & 0.9234 & 0.9217 & 0.9989 & 0.9943 & 0.9966 & 0.9972 & 1.0000 & 0.6388 & 0.9268 & 0.9940 \vspace{1mm}\\
0.9946 & 0.9945 & 0.9953 & 0.9957 & 0.9954 & 0.6785 & 0.6772 & 0.6783 & 0.6783 & 0.6773 & 0.6381 & 0.6387 & 0.6402 & 0.6391 & 0.6388 & 1.0000 & 0.6784 & 0.6384 \vspace{1mm}\\
0.6801 & 0.6822 & 0.6814 & 0.6808 & 0.6806 & 0.9961 & 0.9980 & 0.9958 & 0.9963 & 0.9945 & 0.9273 & 0.9321 & 0.9300 & 0.9288 & 0.9268 & 0.6784 & 1.0000 & 0.9291 \vspace{1mm}\\
0.6399 & 0.6419 & 0.6412 & 0.6406 & 0.6404 & 0.9255 & 0.9273 & 0.9252 & 0.9257 & 0.9240 & 0.9951 & 0.9968 & 0.9952 & 0.9968 & 0.9940 & 0.6384 & 0.9291 & 1.0000 \vspace{1mm}\\
\end{array}
\right)$
}
\]

In this case, $\overline{R}^T$ becomes an equivalence fuzzy relation. Similarly, we can utilize $\overline{R}^T_\lambda$ to classify the fault samples and checked samples. To compare this method with our method, we set $\lambda$ to 0.91, 0.92, 0.93, 0.94, 0.95, 0.96, 0.97, 0.98, 0.99, 1.00. The clustering results are listed in Table 8.
$$\mbox{\bf{\small Table\ 8\ Clustering with $\overline{R}^T$ at different $\lambda$}}\vspace{1mm}$$
  \begin{center}
  \begin{tabular}{cc}
    \toprule[1pt]
    $\lambda$ &  Clustering of the fault samples and check samples \\
    \midrule
   0.91 & $\{1, 2, 3, 4, 5, 16\}, \{6,7, 8, 9, 10,11, 12,  13, 14, 15, 17,18\}$\vspace{1mm}\\
   0.92 & $\{1, 2, 3, 4, 5, 16\}, \{6,7, 8, 9, 10,11, 12,  13, 14, 15, 17,18\}$\vspace{1mm}\\
   0.93 & $\{1, 2, 3, 4, 5, 16\}, \{6,7, 8, 9, 10,12, 13, 17\}, \{11, 12, 13, 14, 15, 17,18\}$\vspace{1mm}\\
    0.94 & $\{1, 2, 3, 4, 5, 16\}, \{6,7, 8, 9, 10,17\}, \{11, 12, 13, 14, 15,18\}$\vspace{1mm}\\
    0.95 & $\{1, 2, 3, 4, 5, 16\}, \{6,7, 8, 9, 10,17\}, \{11, 12, 13, 14, 15,18\}$\vspace{1mm}\\
    0.96 & $\{1, 2, 3, 4, 5, 16\}, \{6,7, 8, 9, 10,17\}, \{11, 12, 13, 14, 15, 18\}$\vspace{1mm}\\
    0.97 & $\{1, 2, 3, 4, 5, 16\}, \{6,7, 8, 9, 10,17\}, \{11, 12, 13, 14, 15,18\}$\vspace{1mm}\\
   0.98 & $\{1, 2, 3, 4, 5, 16\}, \{6,7, 8, 9, 10,17\}, \{11, 12, 13, 14, 15, 18\}$\vspace{1mm}\\
    0.99 & $\{1, 2, 3, 4, 5, 16\}, \{6,7, 8, 9, 10,17\}, \{11, 12, 13, 14, 15, 18\}$\vspace{1mm}\\
    1.00&\{1\},\{2\},\{3\},\{4\},\{5\},\{6\},\{7\},\{8\},\{9\},\{10\},\{11\},\{12\},\{13\},\{14\},\{15\},\{16\},\{17\},\{18\}\vspace{1mm}\\
    \bottomrule[1pt]
  \end{tabular}
\end{center}

Obviously, the clustering results with $\bar{R}^T$ are identical to those obtained using $R$ for $\lambda\in \{0.93,0.94,0.95,$ $0.96,0.97,0.98, 0.99, 1.00\}$. However, samples 17 and 18 cannot be classified by the classical fuzzy clustering method when $\lambda$ ranges from 0.9098 to 0.9288.

Unlike the classical fuzzy clustering method, we directly utilize fuzzy relations to cluster objects. The advantages of our approach can be summarized as follows:
\begin{itemize}
  \item Preserving the features of original data. Our method directly employs the fuzzy relation $R$ for clustering in Example 5.2. In contrast,  classical fuzzy clustering method uses the transitive closure $\overline{R}^T$  as a substitute for $R$. This substitution inevitably introduces distortions to the fuzzy relations inherent in the original data. Following the method in  \cite{He}, we define the distortion degree of two fuzzy relations $R$ and $R'$ on the finite universe $U$\vspace{1mm} as $\delta(R,R')=\sqrt{\mathop{\sum}\limits_{i=1}^n \mathop{\sum}\limits_{j=1}^n \left(R(i,j)-R'(i,j) \right)^2}$. Consequently, we obtain $\delta(R, \overline{R}^T)=0.5786$.
  \item Lower data noise sensitivity. Let us investigate the case when the checked date 17 is corrupted by noise as (0.027, 0.063, 0.760, 0.262, 0.117). In this case, our method can still diagnose it as unbalance when $\lambda=0.9285$. However, the classical fuzzy clustering method fails to make this diagnosis. This demonstrates that our method exhibits moderate noise robustness.
  \item Lower parameter threshold sensitivity. Our method successfully classifies the checked data when $\lambda$ varies  from 0.9098 to 0.9914. However, The classical fuzzy clustering method only works when $\lambda$ varies  from 0.9289 to 0.9914. This means our method achieves a lower $\lambda$ threshold compared to the classical clustering method.
  \item Storage efficiency. In our method, the space complexity of Algorithm 1 is $n^2+1$. Moreover, it needs to update $n^2$ elements to compute $R_\lambda$. Thus, the space complexity of our method is $n^2+1$. However, classical fuzzy clustering method requires $3n^2$\cite{Naessens}. The space complexity of these methods are compared in Table 9.
   \end{itemize}
$$\mbox{\bf{\small Table\ 9 \ Comparing the space complexity}}$$
\begin{center}
 \tabcolsep 0.05in
\begin{tabular}{cc}
 \toprule[1pt]
Fuzzy clustering methods &Space complexity\vspace{1mm}\\
\midrule[0.75pt]
Classical fuzzy clustering method &$3n^2$\vspace{1.5mm}\\
Our method&$n^2+1$\\
   \bottomrule[1pt]
\end{tabular}
\end{center}\vspace{1.5mm}

\qquad However, there exist some limitations in our method as shown in the following.
\begin{itemize}
\item  Deficiency of a theoretical guide on selecting the degree of transitivity. It is clear that valid clustering cannot be achieved with our method when $\alpha_{(I,T)}(R)$ is low. As demonstrated in Example 3.12, where $\alpha_{(I,T)}(R)=0.53$, no $\lambda\in[0,1]$ exists that enables valid clustering of 12 objects. This suggests that a high degree of transitivity is required for a fuzzy relation to cluster objects effectively. However, our method lacks theoretical criteria to determine an appropriate $\alpha_{(I,T)}(R)$.
\item Lack of a theoretical guide on selecting parameter thresholds. Let $T_{\rm M}$ be used to define the transitivity of $R$. It is well-known that $R_\lambda$ is an equivalence relation for any $\lambda\in [0,1]$ when $R$ is an equivalence fuzzy relation. However, this may not hold when other t-norms are chosen to describe the transitivity of $R$. Consequently, objects may fail to be partitioned into cohesive classifications for certain parameter thresholds. Our method currently lacks a systematic approach to selecting these thresholds, which could be a focus of future research.
\end{itemize}
\subsection{Experiments}
\qquad To further validate the advantages of our proposed method, we selected three datasets from the UCI Machine Learning Repository (http://archive.ics.uci.edu) for experimental evaluation. First, 2111 records were selected from the Estimation of Obesity Levels Based on Eating Habits and Physical Condition dataset. Using the cosine similarity method, we obtained the fuzzy relation between the 2111 records. Obviously, this fuzzy relation is not $T_{\rm L}$-transitive. With Algorithm 1, we determined that its degree of transitivity is 0.8739. To compare the classical fuzzy clustering method with our method, its $T_{\rm P}$-transitive closure is then calculated via Algorithm 1 in Ref.\cite{Naessens}. Subsequently, the performance of the classical fuzzy clustering method was compared with our proposed method in Table 10.
$$\mbox{\bf{\small Table\ 10 \ Performance comparison of two methods on Obesity Levels}}\vspace{-4mm}$$
 $$\qquad\qquad\ \mbox{\bf{\small Based On Eating Habitsand Physical Condition dataset}}\vspace{-3mm}$$
\begin{center}
 \tabcolsep 0.05in
\begin{tabular}{cccc}
 \toprule[1pt]
Clustering methods &Running times&Range of $\lambda$&Distortion degree\\
\midrule[0.75pt]
Classical fuzzy clustering &$43.1760''$&[0.9906, 0.9907]&4.3413\vspace{1.5mm}\\
&&$[0.8500,0.8540]\cup[0.8560,0.8589]$&\\
Our method&$15.0057''$&$\cup[0.8627,0.8629]\cup\{0.8659\}$&0\\
&&$\{0.8667\}\cup[0.8710,0.8711]$&\\
   \bottomrule[1pt]
\end{tabular}
\end{center}\vspace{1.5mm}

As demonstrated in Example 5.2, our method requires some time to compute the transitive degree of $R$ by Algorithm 1, whereas classical fuzzy clustering methods spend significantly more time on computing the transitive closure.  Considering the complexities of computation and space (see the discussion in Section 3 and Table 9), Table 10 shows that our method requires less time for clustering compared to classical fuzzy clustering methods.

As mentioned earlier, our method exhibits lower parameter threshold sensitivity. This indicates that that our method successfully clustering samples data when $\lambda\in [0.8500,0.8540]\cup[0.8560,0.8589]\cup[0.8627,0.8629]\cup\{0.8659,0.8667\}\cup[0.8710,0.8711]$.  In contrast, classical fuzzy clustering methods only remain effective within a narrower range of $\lambda$. Moreover, we observe higher distortion in $R$ from Table 10. Consequently, we prefer to cluster the sample data based on the fuzzy relation $R$ in order to enhance the interpretability of clustering from the original data.

We selected the Mice Protein Expression dataset as the second benchmark to evaluate our method's performance. Using $T_{\rm L}$ to measure the transitivity of the fuzzy relation constructed from this dataset, we obtained a degree of transitivity $\alpha_{(I_{\rm GG},T_{\rm L})}(R)=0.8897$. The corresponding performance results are summarized in Table 11.
$$\mbox{\bf{\small Table\ 11 \ Performance comparison of two methods on Mice Protein Expression dataset}}$$
\begin{center}
 \tabcolsep 0.05in
\begin{tabular}{cccc}
 \toprule[1pt]
Clustering methods &Running times&Range of $\lambda$&Distortion degree\\
\midrule[0.75pt]
Classical fuzzy clustering method &$0.2809''$&[0.9745, 0.9753]&13.6000\vspace{1.5mm}\\
&&$[0.9555,0.9587]\cup$&\\
Our method&$0.1786''$&$[0.9622,0.9644]\cup$&0\\
&&[0.9611,0.9688]&\\
   \bottomrule[1pt]
\end{tabular}
\end{center}\vspace{1mm}

Our method's performance on the Statlog (German Credit Data) dataset is likewise summarized in Table 12.
$$\mbox{\bf{\small Table\ 12 \ Performance comparison of two methods on Statlog (German Credit Data) dataset}}$$
\begin{center}
 \tabcolsep 0.05in
\begin{tabular}{cccc}
 \toprule[1pt]
Clustering methods &Running times&Range of $\lambda$&Distortion degree\\
\midrule[0.75pt]
Classical fuzzy clustering method &$1.7625''$&[0.9985, 0.9992]&1.0012\vspace{1.5mm}\\
Our method&$0.8093''$&[0.9939,0.9983]&0\\
   \bottomrule[1pt]
\end{tabular}
\end{center}\vspace{1mm}

Finally, as mentioned in Section 5.2, we emphasize that our method necessitates a high degree of transitivity, which is demonstrated using the 1599 red wine samples from the Wine Quality dataset. With the similar method aforementioned, we can calculate the degree of transitivity $\alpha_{(I_{\rm GG},T_{\rm P})}(R)=0.5987$. However, we cannot find a $\lambda$ such that the 1599 red wine samples are validly clustered. This implies the necessity of determining a suitable threshold for the transitive degree of a fuzzy relation. Investigating how varying degrees of transitivity influence clustering accuracy could represent a future research direction.
\section{Conclusions}
   \qquad  To enhance the application of non-transitive fuzzy relation, this work has investigated the  degrees
of transitivity for a fuzzy relation using several well-known fuzzy implications. Specifically, the conclusions are summarized as follows:

  (1) Given the necessary and sufficient conditions for a fuzzy relation $R$ to be $\varepsilon$-transitive under several well-known fuzzy implications which include $(S,N)$-,
R-, QL-, $g$-implication, $T$-power implication, probabilistic implication
and probabilistic S-implication.

(2)  Investigated the properties of $\alpha_{S_I}(R)$ and its relationship with $\alpha_{(I,T)}(R)$ with several well-known fuzzy implications.

  (3)  Presented an algorithm to calculate the degree of transitivity for a fuzzy relation.

 (4) Characterized the aggregation function $G_n$ which preserves the $\varepsilon$-$T$-transitivity when a series
of $\varepsilon$-$T$-transitive fuzzy relations are aggregated by it.

  (5) Applied the $\varepsilon$-$T$-transitive fuzzy relation to approximate reasoning and clustering.

 (6) Compared our proposed clustering method with the classical fuzzy clustering method.

 (7) Demonstrated the effectiveness of our method using three datasets.

These results provide deeper insights into the degree of transitivity for a fuzzy relation. Future research will explore applications of  $\varepsilon$-transitive fuzzy relation in real-life decision-making.
\section{Acknowledgements}\qquad  This work was funded by the National Natural
Science Foundation of China (Grant No. 61673352) and the Natural Science Basic Research Program of Shaanxi(Grant No. 2023-JC-YB-008).	

\end{document}